\documentclass[10pt]{article}

\usepackage{mypackages}
\usepackage{mycommands}
\usepackage{plot_options}

\title{High-order parallel-in-time method for the monodomain equation in cardiac electrophysiology\thanks{This work was supported by the European High-Performance Computing Joint Undertaking EuroHPC under grant agreements No 955495 (MICROCARD) and No 955701 (TIME-X), co-funded by the Horizon 2020 programme of the European Union (EU) and the Swiss State Secretariat for Education, Research and Innovation.}}

\author[1]{Giacomo Rosilho de Souza\thanks{Corresponding author: \email{giacomo.rosilhodesouza@usi.ch}}}
\author[2]{Simone Pezzuto\thanks{Contributing author: \email{simone.pezzuto@unitn.it}}}
\author[3,4]{Rolf Krause\thanks{Contributing author: \email{rolf.krause@fernhuni.ch}}}

\affil[1]{Euler Institute, Universit\`a della Svizzera italiana, via G. Buffi 13, Lugano, 6900, Switzerland}
\affil[2]{Laboratory for Mathematics in Biology and Medicine, Department of Mathematics, Universit\`a di Trento, via Sommarive 14, Trento, 38123, Italy}
\affil[3]{Faculty of Mathematics and Informatics, FernUni, Schinerstrasse 18, Brig, 3900, Switzerland}
\affil[4]{AMCS, CSE Department, KAUST, Thuwal-Jeddah, Saudi Arabia}

\begin{document}

\maketitle

\begin{abstract}
Simulation of the monodomain equation, crucial for modeling the heart's electrical activity, faces scalability limits when traditional numerical methods only parallelize in space. 
To optimize the use of large multi-processor computers by distributing the computational load more effectively, time parallelization is essential. 
We introduce a high-order parallel-in-time method addressing the substantial computational challenges posed by the stiff, multiscale, and nonlinear nature of cardiac dynamics. 
Our method combines the semi-implicit and exponential spectral deferred correction methods, yielding a hybrid method that is extended to parallel-in-time employing the PFASST framework.
We thoroughly evaluate the stability, accuracy, and robustness of the proposed parallel-in-time method through extensive numerical experiments, using practical ionic models such as the ten-Tusscher--Panfilov. 
The results underscore the method's potential to significantly enhance real-time and high-fidelity simulations in biomedical research and clinical applications.
\end{abstract}



\section{Introduction}\label{sec:intro}

The heart's ability to pump blood hinges on the propagation of electrical signals across cardiac cells and tissues. Central to modeling these complex dynamics is the monodomain equation, which characterizes the propagation of the electric potential and is modulated by ionic models. The latter, describe the interaction of ionic currents and cellular membranes, encompassing a broad range of physiological and pathological cardiac behaviors.

Simulating the monodomain equation is computationally demanding due to its multiscale, nonlinear nature and the high spatio-temporal resolutions required for capturing cardiac dynamics. In traditional numerical methods, optimized for spatial parallel scalability, time integration can be a bottleneck. Indeed, to increase accuracy, both spatial and time resolutions are improved. Although increasing the number of spatial processors keeps the cost per time step constant, it doesn't address the need for additional time steps, and the number of processors is limited by communication bottlenecks. As a result, these schemes have limited scalability, impeding real-time simulations.
To mitigate these limitations, we introduce a high-order parallel-in-time (PinT) method specifically for the monodomain equation. This method enhances scalability by allowing parallelization in the time direction, especially beneficial when spatial scalability is maxed out.

The literature on PinT methods for the monodomain equation is very scarce. In \cite{BedBelHaeGreMouBouBis16}, the authors solve the monodomain equation in parallel-in-time employing the parareal framework \cite{LMT01}. In \cite{BenMinKra21}, the PFASST \cite{EmM12} and Space-Time Multigrid \cite{GanNeu16} methods are compared on a reaction-diffusion problem, which turns out to be a special case of the monodomain equation.
In both cases the scalar nonstiff FitzHugh--Nagumo (FHN) \cite{Fit55,NagAriYos62} ionic model is employed, which is overly simplistic to be employed in practical simulations. Due to stability and efficiency reasons, the methods in \cite{BedBelHaeGreMouBouBis16,BenMinKra21} cannot be used with realistic highly nonlinear and multi-scale models, such as the popular ten-Tusscher--Panfilov (TTP) \cite{tenPan06} or the Courtemance--Ramirez--Nattel (CRN) \cite{CRN98} models, for instance.

After spatial discretization, the monodomain equation can be written as
\begin{subequations}\label{eq:ode_problem}
	\begin{equation+}\label{eq:def_ode}
		y'=f_I(t,y)+f_E(t,y)+f_e(t,y),\qquad y(0)=y_0,
	\end{equation+}
	where
	\begin{equation+}\label{eq:def_fe}
		f_e(t,y)=\Lambda(y)(y-y_\infty(y)).
	\end{equation+}
\end{subequations}
In \cref{eq:ode_problem}, $f_I$ is a stiff linear term representing diffusion of the electric potential, $f_e$ is very stiff and drives the dynamics of the ionic model's gating variables, and $f_E$ is a nonstiff but nonlinear and expensive to evaluate term for the remaining ionic model's variables. The stiff matrix $\Lambda(y)$ is diagonal and $y_\infty(y)$ is a nonlinear term for the gating variables' equilibrium.

The method proposed here leverages the PFASST framework, which integrates the parareal \cite{LMT01}, spectral deferred correction (SDC) \cite{Dutt2000a}, and full approximation scheme (FAS) \cite{TroOosSchBraOswStu07} methods. The parareal algorithm, an iterative PinT method, starts with a coarse solver that progresses through the entire time domain serially; quickly but inaccurately. Then, in parallel, fine solvers refine the results in smaller sub-intervals. 
Enhancements to parareal incorporate the SDC method for both fine and coarse solvers. SDC is an iterative method for solving the monolithic nonlinear system in implicit Runge--Kutta (RK) methods. It works as a preconditioned fixed point iteration and is flexible: able to adjust its preconditioner for specific problems, apply to any RK method, and achieve arbitrary accuracy. When employed in the parareal framework, SDC methods of high and low order serve as fine and coarse solvers, respectively.
Furthermore, SDC synergizes with parareal by blending their iterations and improving parallel efficiency.
Combined with FAS to ensure precise information transfer between levels, this gives rise to the multi-level PinT PFASST method \cite{EmM12,EmmMin14,BolMosSpe17}.

Based on this framework, new PinT methods are designed by replacing SDC with one of its variants. For instance, the semi-implicit SDC (SI-SDC) \cite{Min03}, for problems \cref{eq:ode_problem} with $f_e\equiv 0$. Or the exponential SDC (ESDC) \cite{Buv20}, for \cref{eq:ode_problem} with $f_I\equiv f_E\equiv 0$. In this paper, we propose a new hybrid SDC (HSDC) method that combines SI-SDC, for $f_I$ and $f_E$, with ESDC, for $f_e$. A PinT method based on HSDC is seamlessly obtained using the PFASST methodology. For simplicity, we will denote the method HSDC in both its serial and parallel versions. 

Although the HSDC method is designed for the monodomain equation, it can be applied to any problem structured as \cref{eq:ode_problem}. Like SI-SDC and ESDC, it has an arbitrary order of convergence. This is uncommon for the monodomain equation even in the serial context; indeed, for efficiency and stability, most methods employ splitting strategies and are low order. Further, due to exponential integration, the method is stable even for severely stiff ionic models, such as TTP. Finally, the method is computationally cost-effective, as its preconditioner relies on the efficient Rush--Larsen method \cite{RuL78,LinGerJahLoeWeiWie23} popular in cardiac electrophysiology; that is, a splitting method based on implicit, explicit, and exponential Euler steps. The HSDC method proposed here has been implemented in \texttt{pySDC} \cite{Spe19}, a Python package for prototyping and test SDC-like methods.

The paper is structured as follows. \Cref{sec:monodomain} reviews the monodomain equation and \cref{sec:pint_monodomain} describes the development of the PinT HSDC method. In \cref{sec:stab} we perform a numerical stability analysis of the HSDC method.
In \cref{sec:numexp}, we present a comprehensive set of numerical experiments that demonstrate the robustness and effectiveness of HSDC. Conclusions are drawn in \cref{sec:conclusion}.

\section{The monodomain equation for cardiac electrophysiology}\label{sec:monodomain}
The monodomain equation is a fundamental tool in cardiac electrophysiology, describing how electrical signals propagate through heart tissue. It models the interplay of ion channels, cellular membranes, tissue conductivity and electric potential propagation. Therefore, it provides insights into phenomena ranging from normal cardiac rhythm to pathological conditions such as arrhythmia.

\subsection{The monodomain model}
The monodomain equation is given by the following reaction-diffusion system of equations
\begin{subequations}\label{eq:monodomain}
	\begin{align+}
		C_m \frac{\partial V_m}{\partial t} &= \chi^{-1} \nabla \cdot \bigl(\sigma \nabla V_m\bigr) + I_\text{stim}(t,\bm x) - I_\text{ion}(V_m,{w}_a,{w}_g), \qquad && \text{in $\Omega\times (0,T]$,} \label{eq:monodomain_pde} \\		
		\frac{\partial{w}_a}{\partial t} &= h_a(V_m,{w}_a,{w}_g), && \text{in $\Omega\times (0,T]$,} \label{eq:monodomain_a} \\
		\frac{\partial{w}_g}{\partial t} &= h_g(V_m,{w}_g) = {\alpha}(V_m)\bigl(\mathbf{1}-{w}_g\bigr) - {\beta}(V_m) {w}_g, && \text{in $\Omega\times (0,T]$,} \label{eq:monodomain_g} \\
		-\sigma\nabla V_m\cdot \bm{n} &= 0, && \text{on }\partial \Omega\times (0,T], \\
		V_m(0,\bm x) &= V_{m,0}, \quad {w}_a(0,\bm x) = {w}_{a,0}, \quad {w}_g(0,\bm x) = {w}_{g,0}, && \text{in $\Omega$}
	\end{align+}
\end{subequations}
where $T>0$ and $\Omega\subset\Rb^d$ is a domain. 
The model unknowns are the the transmembrane potential $V_m:\Rb^d\times\Rb\rightarrow\Rb$, in \si{\milli\volt}, the vector of auxiliary ionic variables $w_a:\Rb^d\times\Rb\rightarrow\Rb^{m_1}$, and the vector of ionic gating variables $w_g:\Rb^d\times\Rb\rightarrow\Rb^{m_2}$. 
The values for the isotropic conductivity coefficient $\sigma$, the surface area-to-volume ratio $\chi$ and the cellular membrane electrical capacitance $C_m$ are given in \cref{sec:numexp}.

The $I_\text{stim}$ term is used to stimulate the tissue and initiate an action potential's propagation.
The ionic current $I_\text{ion}$ models the transmembrane currents and, like $h_a$, $\alpha$, $\beta$ and $m_1$, $m_2$, its exact definition depends on the chosen ionic model. Several models exists in the literature, and the choice depends on the modeled species, type of muscle and pathological conditions, for instance. Nevertheless, in most cases the structure is as in \cref{eq:monodomain}, hence with $\mathbf{1}\in\Rb^{m_2}$ a vector of ones, $\alpha(V_m)$, $\beta(V_m)$ diagonal matrices and $h_a$, $I_\text{ion}$ general nonlinear terms. Therefore,
\begin{equation}\label{eq:monodomain_g_lambda}
	\frac{\partial w_g}{\partial t} = \Lambda_g(V_m)\bigl({w}_g - {w}_{g,\infty}(V_m) \bigr), 
\end{equation}
where ${\Lambda}_g(V_m) = -({\alpha}(V_m) + {\beta}(V_m))$ is a diagonal matrix and ${w}_{g,\infty}(V_m) =({\alpha}(V_m) + {\beta}(V_m))^{-1}{\alpha}(V_m)\mathbf{1}$ is the steady state value of ${w}_g$ (for fixed $V_m$).

\subsection{Semi-discrete equation}\label{sec:monodomain_space_disc}
For the purpose of this paper we consider rectangular domains $\Omega$ and discretize the monodomain equation \cref{eq:monodomain} using a fourth-order finite difference scheme, yielding the semi-discrete system
\begin{subequations}\label{eq:sd_sys}
	\begin{align+}\label{eq:sd_sys_a}
		C_m\frac{\dif \bm V_m}{\dif t} &= \chi^{-1}\mathbf{A}\bm{V}_m +\bm{I}_\text{stim}(t) - \bm{I}_\text{ion}(\bm{V}_m,\bm{w}_a,\bm{w}_g),\\ \label{eq:sd_sys_b}
		\frac{\dif \bm{w}_a}{\dif t} &= \bm{h}_a(\bm{V}_m,\bm{w}_a,\bm{w}_g),\\ \label{eq:sd_sys_c}
		\frac{\dif \bm{w}_g}{\dif t} &= \bm{\Lambda}_g(\bm{V}_m)(\bm{w}_g-\bm{w}_{g,\infty}(\bm{V}_m)),
	\end{align+}
\end{subequations}
plus the initial conditions, where $\bm{V}_m:\Rb\rightarrow\Rb^{N_{dof}}$, $\bm{w}_a:\Rb\rightarrow\Rb^{m_1 N_{dof}}$, $\bm{w}_g:\Rb\rightarrow\Rb^{m_2 N_{dof}}$ and $N_{dof}$ is the number of mesh degrees of freedom.
The ionic model's terms $\bm{h}_a,\bm{\Lambda}_g,\bm{w}_{g,\infty},\bm I_\text{ion}$ and the stimulus $\bm{I}_{\text{stim}}$ are evaluated node-wise.

The semi-discrete equation \cref{eq:sd_sys} can be written as \cref{eq:ode_problem}
for $y=(\bm{V}_m,\bm{w}_a,\bm{w}_g)\in\Rb^N$, with $N=(1+m_1+m_2)N_{dof}$, by defining
\begin{subequations}\label{eq:def_monodomain_rhs}
	\begin{gather+}\label{eq:def_f_IE}
		f_I(t,y)=\begin{pmatrix}
			(\chi C_m)^{-1} \mathbf{A}\bm{V}_m \\ \bm{0} \\ \bm{0}
		\end{pmatrix}, \qquad
		f_E(t,y)=\begin{pmatrix}
			C_m^{-1}(\bm{I}_{stim}(t)-\bm{I}_{ion}(\bm{V}_m,\bm{w}_a,\bm{w}_g)) \\  \bm{h}_a(\bm{V}_m,\bm{w}_a,\bm{w}_g) \\ \bm{0}\end{pmatrix}, 
	\end{gather+}
	and $f_e(t,y)=\Lambda(y)(y-y_\infty(y))$, where
	\begin{equation+}\label{eq:def_Lambda}
		\Lambda(y)=\begin{pmatrix}
			\bm{0} & \bm{0} & \bm{0} \\ \bm{0} & \bm{0} & \bm{0} \\ \bm{0} & \bm{0} & \bm{\Lambda}_g(\bm{V}_m)
		\end{pmatrix}, \quad
		y_\infty(y) = \begin{pmatrix}
			\bm{0} \\ \bm{0}\\ \bm{w}_{g,\infty}(\bm{V}_m)
		\end{pmatrix}
	\end{equation+}
\end{subequations}
and $\bm{\Lambda}_g(\bm{V}_m)$ is a diagonal matrix.

Concerning the numerical integration of \cref{eq:sd_sys} written in the form of \cref{eq:ode_problem}, note that $f_I$ is stiff but linear and $\mathbf{A}$ is symmetric, hence it can be integrated relatively efficiently by an implicit method. 
The $f_E$ term houses the most intricate components of the ionic model, hence it is expensive. Yet, it is nonstiff and can therefore be integrated explicitly. Finally, $f_e$ is very stiff but due to the diagonal form of $\Lambda(y)$ it is easily integrated exponentially. These considerations motivate the Rush--Larsen method, proposed in \cite{RuL78}, where a splitting strategy is used and $f_I$, $f_E$, and $f_e$ are integrated by the implicit, explicit and exponential Euler methods, respectively.
This approach is still very popular nowadays, due to its robustness, efficiency and simple implementation \cite{LinGerJahLoeWeiWie23}. In \cref{sec:prec_picard} we propose a preconditioner based on this approach.

\section{Parallel-in-time method for the monodomain equation}\label{sec:pint_monodomain}
In this section, we present the PinT method for the monodomain equation defined by \cref{eq:ode_problem,eq:def_monodomain_rhs}. 
To design such method, we first formulate a new SDC method where the three terms $f_I$, $f_E$, and $f_e$ are integrated with appropriate schemes.
Then we derive its PinT extension using the PFASST framework, to do so we follow the established methodology \cite{BolMosSpe18,BolMosSpe17,SanSchSpe22,SchSpe20,Spe19} employing preconditioned Picard iterations and the multigrid perspective.

\subsection{Hybrid spectral deferred correction method}
Here we derive the HSDC method for the monodomain equation. This method is an hybrid between ESDC \cite{Buv20} and SI-SDC \cite{Min03}, where the exponential term $f_e$ is integrated with ESDC and the implicit and explicit terms $f_I,f_E$ are integrated with the SI-SDC method. The construction of such method follows three basic steps. First, the monodomain equation is written in a form more compatible with exponential Runge--Kutta (ERK) methods \cite{HoO10,HocOst2005}. Then ERK methods of collocation type \cite{HocOst2005} are re-derived, we do that for completeness and mostly to obtain a formulation 
closer to collocation RK methods and therefore easily implementable in available software for PinT methods, such as \texttt{pySDC} \cite{Spe19}. Finally, we recall the preconditioned Picard iteration for solving the nonlinear system arising in each step of the ERK method and propose a preconditioner tailored to the monodomain equation.

\subsubsection{Monodomain equation reformulation}
Here we rewrite the monodomain equation in a ERK-friendly formulation.

Let $\Delta t>0$ be the step size and $t_n=n\Delta t$.
In the time interval $t\in[t_n,t_{n+1}]$, equation \cref{eq:ode_problem} with initial value $y_n$ can be written as
\begin{equation}\label{eq:def_ode_esdc}
	y'=\Ln (y-y_n) + \Gn(y),\qquad y(t_n)=y_n,
\end{equation}
where
\begin{equation}\label{eq:def_g}
	\Gn(y)=f_I(t,y)+f_E(t,y)+f_e(t,y)+ \Ln(y_n-y)
\end{equation}
and, for simplicity, $y_n$ and the time variable $t$ are dropped from the arguments of $g$.

In the next section we define the ERK methods for equation \cref{eq:ode_problem} written in the form \cref{eq:def_ode_esdc}, where the matrix exponential of $\Ln$ will be needed. Concerning accuracy, $\Ln$ could be replaced by any matrix; however, for stability reasons, this matrix should capture the stiffness of the problem. In our case we have chosen $\Ln$, so that it captures the stiffness of the exponential term $f_e$ evaluated at values $y$ sufficiently close to $y_n$. In contrast, the stiffness of the $f_I$ term will be resolved implicitly. This choice is particularly efficient since due to the diagonal form of $\Ln$, leading to cheap matrix exponential. Additionally, $f_I$ is linear, hence its implicit integration is as well relatively inexpensive.

\subsubsection{Exponential Runge--Kutta methods}\label{sec:erk}
Here we derive ERK methods of collocation type for \cref{eq:def_ode_esdc}, but choose a formulation resembling collocation RK methods.

Applying the variation of constants formula to \cref{eq:def_ode_esdc} we obtain\footnote{Notice that, generally, there is an exponential in front of $y_n$ \cite{HoO10,HocOst2005}. Here, we prefer to put the exponential inside the integral and take out the identity operator. This approach leads to a formulation of ERK closer to standard RK methods.}
\begin{equation}
	y(t)=y_n+\int_{t_n}^t e^{(t-s)\Ln}\Gn(y(s)) \dif s,
\end{equation}
and performing the change of variables $t=t_n+r\Delta t$ yields
\begin{equation}\label{eq:var_cte}
	y(t_n+r\Delta t)=y_n+\Delta t\int_{0}^r e^{(r-s)\Delta t\Ln}\Gn(y(t_n+s\Delta t)) \dif s.
\end{equation}
To discretize \cref{eq:var_cte} and obtain exponential integrators, $\Gn(y(t_n+s\Delta t))$ is approximated by an interpolating polynomial.
Let $0<c_1<\ldots<c_M=1$ be Radau IIA collocation nodes and, for notational purposes, we also define $c_0=0$. We denote $\by_n=(y_{n,0},\ldots,y_{n,M})\in\Rb^{N(M+1)}$, where the $y_{n,j}\in\Rb^N$ represent approximations to $y(t_{n,j})$, with $t_{n,j}=t_n+c_j\Delta t$.
Let
\begin{equation}\label{eq:piN}
	\begin{aligned}
		p_{\Gn}(\by_n)(s)&= \sum_{j=1}^M \Gn(y_{n,j})\ell_j(s),&\text{with}&& \ell_j(s)&=\prod_{\substack{k=1\\k\neq j}}^M \frac{s-c_k}{c_{j}-c_{k}},
	\end{aligned}
\end{equation}
be a polynomial interpolating $\Gn(y_{n,j})$, $j=1,\ldots,M$. Note that $y_{n,0}$ and $c_0$ are excluded from the definitions of $p_{\Gn}(\by_n)$ and the Lagrange basis functions $\ell_j$.
Replacing $\Gn(y(t_n+s\Delta t))$ with $p_{\Gn}(\by_n)(s)$ in \cref{eq:var_cte} and evaluating at $r=c_i$ yields the system of equations for $y_{n,i}$:
\begin{equation}\label{eq:sys_erk}
	y_{n,i} = y_n+\Delta t\int_{0}^{c_i} e^{(c_i-s)\Delta t\Ln}p_{\Gn}(\by_n)(s) \dif s, \qquad i=0,\ldots,M.
\end{equation}
Let 
\begin{equation}\label{eq:def_aij}
	a_{ij}(z)=\int_0^{c_i}e^{(c_i-s)z}\ell_j(s)\dif s,
\end{equation}
then \cref{eq:sys_erk} is rewritten as
\begin{equation}\label{eq:sys_erk_aij}
	y_{n,i} = y_n+\Delta t\sum_{j=1}^M a_{ij}(\Delta t\Ln)\Gn(y_{n,j}), \qquad i=0,\ldots,M,
\end{equation}
and the approximation of $y(t_{n+1})$ is given by $y_{n+1}=y_{n,M}$, since $c_M=1$.
A stable way of computing the matrix coefficients $a_{ij}(\Delta t \Ln)$ is proposed in \cref{app:comp_aij}.

Remark that the first equation of \cref{eq:sys_erk_aij}, $i=0$, is simply $y_{n,0}=y_n$. Therefore, system \cref{eq:sys_erk_aij} can be written in compact form as a $N(M+1)\times N(M+1)$ nonlinear system
\begin{equation}\label{eq:def_C}
	\bm{C}(\by_n)\coloneqq \left(I -\Delta t A(\Delta t \Lnz )\bG\right)(\by_n) =\mathbf{1}\otimes y_n,
\end{equation}
where $\otimes$ is the Kronecker product,
$\mathbf{1}\in\Rb^{M+1}$ a vector of ones, and
$I\in\mathbb{R}^{N(M+1)\times N(M+1)}$ is the identity matrix. 
The block matrix $A(\Delta t \Lnz )\in\mathbb{R}^{N(M+1)\times N(M+1)}$, with $(M+1)\times (M+1)$ blocks, and the nonlinear term $\bG:\Rb^{N(M+1)}\rightarrow\Rb^{N(M+1)}$, are given by
\begin{equation}\label{eq:def_A_matrix}
	A(z) = \begin{pmatrix}
		0 & 0 & \cdots & 0 \\
		0 & a_{11}(z) & \cdots & a_{1M}(z) \\
		\vdots & \vdots & \ddots & \vdots \\
		0 & a_{M1}(z) & \cdots & a_{MM}(z)
	\end{pmatrix}, \qquad
	\bG(\by_n)=
	\begin{pmatrix}
		\Gn(y_{n,0}) \\ \Gn(y_{n,1}) \\ \vdots \\ \Gn(y_{n,M}))
	\end{pmatrix},
\end{equation}
where in the definition \cref{eq:def_g} of $g$ we also replace $y_n$ with $y_{n,0}$. 
Indeed, the main advantage of including $y_{n,0}$ into $\by_n$ is that due to the first line of \cref{eq:def_C} we can interchange $y_{n,0}$ and $y_n$. In this way the definition of $\bm C$ doesn't directly depend on $y_n$ and can be used for all time steps.

A few comments are warranted regarding the interconnections between \cref{eq:def_C}, collocation methods and exponential Runge--Kutta methods. First, for $z=0$, \cref{eq:def_aij} reduce to the weights of a standard collocation method \cite[Section II.7]{HNW08}.
Here $z=\Delta t\Lnz$ is a nonzero block diagonal matrix \cref{eq:def_Lambda}; however, the first two diagonal blocks of $\Lnz$ are indeed zero matrices. Thus, for the first two components $V_m,\bm{w}_a$ of $y$, \cref{eq:def_C} simplifies a collocation scheme, which motivates our choice of Radau IIA nodes for $c_1,\ldots,c_M$. 
In contrast, the last diagonal block $\Lambda_g(V_m)$ of $\Lnz$ is nonzero, thus for the last component $\bm{w}_g$ of $y$, \cref{eq:def_C} is an exponential Runge--Kutta method of collocation type \cite{HocOst2005}. Recall that $\Lambda_g(V_m)$ is diagonal, therefore its exponentials are inexpensive to compute. We conclude that, for the monodomain equation, method \cref{eq:def_C} is an hybrid between collocation schemes and exponential Runge--Kutta methods.

\subsubsection{Preconditioned Picard iteration}\label{sec:prec_picard}
The monolithic $N(M+1)\times N(M+1)$ nonlinear system \cref{eq:def_C} is expensive to solve. Furthermore, for some initial guesses, the Newton method might struggle to converge.
Spectral deferred corrections are alternative iterative methods for solving \cref{eq:def_C}. 
These methods were introduced in \cite{Dutt2000a}, where they have been derived by defining a correction equation, which is iteratively solved and improves the current solution (see also  \cite{Buv20,BouLayMin03,Min03,HagZho06,Wei15,SpeRupMinEmmKra16,WeiGho18,CauSea19}). More recently, the same methods were derived by solving a preconditioned Picard iteration \cite{BolMosSpe17,BolMosSpe18,Spe19,KreSpe21,SchSpe20,SanSchSpe22}.
For the purpose of this paper, we follow the latter approach, due to its convenient notation for the parallel-in-time extension.

The preconditioned Picard iteration for \cref{eq:def_C} is given by
\begin{equation}\label{eq:prec_picard}
	\bm{P}(\by_n^{k+1})=(\bm{P}-\bm{C})(\by_n^{k})+\mathbf{1}\otimes y_n,
\end{equation}
where $\bm{P}$ is the nonlinear preconditioner operator and as initial guess we choose $\by_n^0=\mathbf{1}\otimes y_n$. In the remainder of this section we define the preconditioner $\bm{P}$ for the monodomain equation. 

From \cref{eq:var_cte}, we have
\begin{align}
	y(t_n+c_i\Delta t)&=y_n+\Delta t\int_{0}^{c_i} e^{(c_i-s)\Delta t\Ln}\Gn(y(t_n+s\Delta t)) \dif s\\
	&=y_n+\Delta t\sum_{j=1}^i e^{(c_i-c_j)\Delta t\Ln}\int_{c_{j-1}}^{c_j} e^{(c_j-s)\Delta t\Ln }\Gn(y(t_n+s\Delta t)) \dif s.
\end{align}
Let $d_i=c_i-c_{i-1}$ for $i=1,\ldots,M$. 
Since $e^{(c_i-c_j)\Delta t\Ln}=1+\mathcal{O}(\Delta t)$, these terms can be dropped from the forthcoming computations without affecting accuracy. Therefore, approximations $\tilde y_{n,i}$ to $y(t_n+c_i\Delta t)$, for $i=0,\ldots,M$, are obtained as
\begin{equation}\label{eq:deriv_P_a}
	\begin{aligned}
		\tilde y_{n,i}&=y_n+\Delta t\sum_{j=1}^i\int_{c_{j-1}}^{c_j} e^{(c_j-s)\Delta t\Ln}(f_I(\tilde y_{n,j})+f_E(\tilde y_{n,j-1})+f_e(\tilde y_{n,j-1})+\Ln(y_n-\tilde y_{n,j-1})) \dif s\\
		&= y_n+\Delta t\sum_{j=1}^i  d_j\varphi_1(d_j\Delta t \Ln)(f_I(\tilde y_{n,j})+f_E(\tilde y_{n,j-1})+f_e(\tilde y_{n,j-1})+\Ln(y_n-\tilde y_{n,j-1})),
	\end{aligned}
\end{equation}
with
\begin{equation}
	\varphi_1(z)=\int_0^1 e^{(1-r)z}\dif r=\frac{e^z-1}{z}.
\end{equation}
The first two diagonal blocks of $\Ln$ are zero and $\varphi_1(0)=1$, therefore, from \cref{eq:def_f_IE},
\begin{equation}
	\varphi_1(d_j\Delta t \Ln)(f_I(\tilde y_{n,j})+f_E(\tilde y_{n,j-1}))=f_I(\tilde y_{n,j})+f_E(\tilde y_{n,j-1})
\end{equation}
and \cref{eq:deriv_P_a} is 
\begin{equation}\label{eq:deriv_P_b}
	\tilde y_{n,i}= y_n+\Delta t\sum_{j=1}^i  d_j(f_I(\tilde y_{n,j})+f_E(\tilde y_{n,j-1})+\varphi_1(d_j\Delta t \Ln)(f_e(\tilde y_{n,j-1})+\Ln(y_n-\tilde y_{n,j-1}))).
\end{equation}
As previously mentioned, neglecting terms $e^{(c_i-c_j)\Ln\Delta t}=1+\mathcal{O}(\Delta t)$ do not impact accuracy; indeed, the introduced error is comparable, asymptotically, to the one of the approximated integrals. The reason for their removal is to improve computational efficiency, since they would require additional $\mathcal{O}(M^2)$ diagonal matrix exponential computations. Furthermore, \cref{eq:deriv_P_b} corresponds to the popular IMEX Rush--Larsen method in cardiac electrophysiology, where the diffusion is integrated implicitly, the stiff gating variables exponentially, and the remaining nonstiff terms explicitly \cite{RuL78,LinGerJahLoeWeiWie23}.

Overall, \cref{eq:deriv_P_b} is significantly cheaper to compute than \cref{eq:sys_erk_aij}, indeed it is a diagonally implicit method where the implicit term is linear.
Therefore, from \cref{eq:deriv_P_b}, we define the preconditioning operator $\bm{P}:\Rb^{N(M+1)}\rightarrow\Rb^{N(M+1)}$ as $\bm{P}(\by_n)=(\bm{P}_{0}(\by_n),\ldots,\bm{P}_{M}(\by_n))$, with
\begin{equation}\label{eq:def_P}
	\begin{aligned}
		\bm{P}_{i}(\by_n) = y_{n,i}-\Delta t\sum_{j=1}^i d_j&\left(f_I(y_{n,j})+ f_E(y_{n,j-1})\right.\\
		&\left.\quad+\varphi_1(d_j\Delta t \Ln)(f_e(y_{n,j-1})+\Ln (y_n-y_{n,j-1})) \right)
	\end{aligned}
\end{equation}
for $i=0,\ldots,M$. Remark that for the $f_I,f_E$ terms this is the standard preconditioner used in SI-SDC methods \cite{Min03}. However, here we also have exponential integration of the $f_e$ term.

With a close look to \cref{eq:prec_picard,eq:def_P}, we notice that $y_{n,0}^{k}=y_n$ for all $k$. Hence, for computational efficiency, all terms in \cref{eq:def_P} depending on $y_{n,0}$ can be neglected, since they appear on both sides of \cref{eq:prec_picard}.

	
	\subsection{Parallel-in-time method}\label{sec:pint}
	From \cref{eq:prec_picard}, the derivation of a PinT method based on the PFASST methodology \cite{EmM12,EmmMin14} is obtained following the standard procedure based on the Full Approximation Scheme (FAS) \cite{BolMosSpe18,BolMosSpe17,SanSchSpe22,SchSpe20,Spe19}, which we recall here for completeness.
	
	A sequence of $P$ consecutive time steps of the ERK problem \cref{eq:def_C} can be written as a $N(M+1)P\times N(M+1)P$ nonlinear system
	\begin{align}\label{eq:erk_sys}
		\begin{pmatrix}
			\bm{C} & & &\\
			-\bm{H} & \bm{C} & & \\
			& \ddots &\ddots & \\
			& & -\bm{H} & \bm{C}
		\end{pmatrix}
		\begin{pmatrix}
			\by_0 \\ \by_1 \\ \vdots \\ \by_{P-1}
		\end{pmatrix}
		= 
		\begin{pmatrix}
			\mathbf{1}\otimes y_0 \\ \bm{0}\\ \vdots \\ \bm{0}
		\end{pmatrix},
	\end{align}
	with $\bm{H}=\bm{N}\otimes I_N$ the matrix taking the last node value at the previous step to be used as initial value in the current step, hence $\bm{N}\in\Rb^{(M+1)\times (M+1)}$ has ones in the last column and zeros elsewhere.
	
	Let $\bz=(\by_0,\ldots,\by_{P-1})\in\Rb^{N(M+1)P}$ and $\bm{C}_P(\bz)=(\bm{C}(\by_0),\ldots,\bm{C}(\by_{P-1}))^\top$, then \cref{eq:erk_sys} can be written as
	\begin{equation}\label{eq:def_D}
		\bm{D}(\bz)\coloneqq (\bm{C}_P-\bm{E}\otimes\bm{H})(\bz)=\bm{b},
	\end{equation}
	with $\bm{E}\in\Rb^{P\times P}$ having ones in the first lower off-diagonal and zeros elsewhere, and $\bm{b}=(\mathbf{1}\otimes y_0,\bm{0},\ldots,\bm{0})\in\Rb^{N(M+1)P}$.
	In order to solve \cref{eq:def_D} with a preconditioned Picard iteration, hence
	\begin{equation}\label{eq:prec_picard_seq}
		\bm{Q}(\bz^{k+1})=(\bm{Q}-\bm{D})(\bz^{k})+\bm{b},
	\end{equation}
	there are two types of preconditioners $\bm{Q}$\footnote{Throughout this section, the global ($P$-steps) version of local ($1$-step) operators is denoted by the subsequent letter; for instance, the global version of $\bm{C}$ is denoted $\bm{D}$, and the one of $\bm{P}$ is $\bm{Q}$.}:
	\begin{equation}\label{eq:def_Q}
		\begin{aligned}		
			\bm{Q}^{seq}(\bz)&=(\bm{P}_P-\bm{E}\otimes\bm{H})(\bz),\\
			\bm{Q}^{par}(\bz)&=\bm{P}_P(\bz),
		\end{aligned}
	\end{equation}
	with $\bm{P}_P(\bz)=(\bm{P}(\by_0),\ldots,\bm{P}(\by_{P-1}))^\top$. Notice as $\bm{Q}^{par}$ can be applied in parallel across the time steps, while $\bm{Q}^{seq}$ propagates sequentially the solution form one step to the next. 
	
	The PFASST parallel-in-time method employs the full approximation scheme (FAS) \cite{TroOosSchBraOswStu07} to solve \cref{eq:prec_picard_seq} on a multilevel hierarchy, employing the sequential preconditioner $\bm{Q}^{seq}$ at the coarsest level and the parallel one $\bm{Q}^{par}$ at any other level. When going from one level to the next, coarsening can happen in both space and time. Space coarsening is usually done by decreasing the mesh resolution and/or the polynomial order (as in $hp$-refinement). Time coarsening is done by keeping the step size $\Delta t$ constant but decreasing the number $M$ of collocation nodes. For the purpose of this paper, we focus on time coarsening only.
	
	Let $M_L<M_{L-1}<\ldots<M_1$ be a sequence of collocation nodes' numbers and $\{c_i^l\}_{i=1}^{M_l}$ the corresponding Radau IIA collocation nodes, for $l=1,\ldots,L$. Let $\bm{D}_l$, and $\bm{Q}_l^{seq},\bm{Q}_l^{par}$ be as in \cref{eq:def_D,eq:def_Q}, respectively, but for $M=M_l$ collocation nodes. Similarly, let $\by_{n}^l=(y_{n,0}^l,\ldots,y_{n,M_l}^l)\in\Rb^{N(M_l+1)}$ and $\bz^l=(\by_0^l,\ldots,\by_{P-1}^l)\in\Rb^{N(M_l+1)P}$.
	
	Let $\bm{I}_l:\Rb^{N(M_l+1)}\times \Rb\rightarrow \Rb^N$ be the Lagrange interpolator at nodes $c_1^l,\ldots,c_{M_l}^l$ defined by
	\begin{equation}
		\bm{I}_l(\by_n^l,s)=\sum_{j=1}^{M_l}\ell_j^l(s)y_{n,j}^l,\qquad \ell_j^l(s)=\prod_{\substack{k=1\\ k\neq j}}^{M_l}\frac{s-c_k^l}{c_j^l-c_k^l}.
	\end{equation}
	The restriction operator $\bm{R}_l^{l+1}:\Rb^{N(M_l+1)}\rightarrow\Rb^{N(M_{l+1}+1)}$, from level $l$ to level $l+1$, is linear and we represent it by a matrix, implicitly defined by
	\begin{equation}
		\bm{R}_l^{l+1}\by_n^l\coloneqq (y_{n,0}^l,\bm{I}_l(\by_n^l,c_1^{l+1}),\ldots,\bm{I}_l(\by_n^l,c_{M_{l+1}}^{l+1}))^\top.
	\end{equation}
	The prolongation operator $\bm{T}_{l+1}^l:\Rb^{N(M_{l+1}+1)}\rightarrow\Rb^{N(M_{l}+1)}$ is defined as
	\begin{equation}
		\bm{T}_{l+1}^l\by_n^{l+1}\coloneqq (y_{n,0}^{l+1},\bm{I}_{l+1}(\by_n^{l+1},c_1^{l}),\ldots,\bm{I}_{l+1}(\by_n^{l+1},c_{M_{l}}^{l}))^\top.
	\end{equation}
	The global restriction operator $\bm{S}_l^{l+1}:\Rb^{N(M_l+1)P}\rightarrow\Rb^{N(M_{l+1}+1)P}$ is simply $\bm{S}_l^{l+1}=I_P\otimes \bm{R}_l^{l+1}$. Similarly, the global prolongation operator is $\bm{U}_{l+1}^l=I_P\otimes \bm{T}_{l+1}^l$.

	In \cref{algo:hsdc} we provide a pseudo code for the multilevel and parallel HSDC method. Therein, the \textsc{Burn-in} method is an initialization procedure described in \cite{EmM12}. We stress that \cref{algo:hsdc} is a rough description of the true computational flow of HSDC, which employs the implementation available in the \texttt{pySDC} \cite{Spe19} library and is slightly different. For instance, in the true code, if all the residuals up to some step $n$ are below the tolerance, then iterations at those steps stop and only the remaining steps continue iterating.
	
	\begin{algorithm}
		\begin{algorithmic}[1]
			\Function{HSDC}{$y_0$, $\Delta t$, $P$, $\{M_l\}_{l=1}^L$, $K$}
			\State $\bz_1^0=$ \Call{Burn-in}{$y_0$, $\Delta t$, $P$, $\{M_l\}_{l=1}^L$} 
			\State $\bm{b}_1 = (\mathbf{1}\otimes y_0,\bm{0},\ldots,\bm{0})$ and $\bm{b}_l=\bm{S}_{l-1}^l\bm{b}_{l-1}$, $l=2,\ldots,L$.
			\State $k=0$, \textbf{converged} = False
			\While{not \textbf{converged} and $k<K$}
			\State $\bz_1^{k+1/2}=\bz_1^{k}$, $\bm{\tau}_1^k=\bm{0}$
			\For{$l=2,\ldots,L-1$}
			\State $\bm{\tau}_l^k =\bm{D}_l(\bm{S}_{l-1}^l\bz_{l-1}^{k+1/2}) -\bm{S}_{l-1}^l\bm{D}_{l-1}(\bz_{l-1}^{k+1/2})+\bm{S}_{l-1}^l\bm{\tau}_{l-1}^k$
			\Comment{Residual correction}
			\State Solve $\bm{Q}_l^{par}(\bz_l^{k+1/2})=(\bm{Q}_l^{par}-\bm{D}_l)(\bm{S}_{l-1}^l\bz_{l-1}^{k+1/2})+\bm{b}_l+\bm{\tau}_l^k$
			\EndFor
			\State $\bm{\tau}_L^k =\bm{D}_L(\bm{S}_{L-1}^L\bz_{L-1}^{k+1/2}) -\bm{S}_{L-1}^L\bm{D}_{L-1}(\bz_{L-1}^{k+1/2})+\bm{S}_{L-1}^L\bm{\tau}_{L-1}^k$
			\State Solve $\bm{Q}_L^{seq}(\bz_L^{k+1})=(\bm{Q}_L^{seq}-\bm{D}_L)(\bm{S}_{L-1}^L\bz_{L-1}^{k+1/2})+\bm{b}_L+\bm{\tau}_L^k$
			\For{$l=L-1,\ldots,1$}
			\State $\bm{c}_{l}^{k+1/2} = \bm{T}_{l+1}^l(\bz_{l+1}^{k+1}-\bm{S}_{l}^{l+1}\bz_l^{k+1/2})$
			\Comment{Coarse grid correction}
			\State Solve $\bm{Q}_l^{par}(\bz_l^{k+1})=(\bm{Q}_l^{par}-\bm{D}_l)(\bz_l^{k+1/2}+\bm{c}_{l}^{k+1/2})+\bm{b}_l+\bm{\tau}_l^k$
			\EndFor	
			\State $k\leftarrow k+1$
			\State $(\bm{r}_0^k,\ldots,\bm{r}_{P-1}^k) \coloneqq \bm{D}_1(\bz^k_1)-\bm{b}_1$
			\Comment{Residuals}
			\State $(\by_0^k,\ldots,\by_{P-1}^k)=\bz^k_1$
			\Comment{Solutions}
			\State \textbf{converged} = True if $\Vert \bm{r}_n^k\Vert<tol\cdot \Vert\by_{n}^k\Vert$ for $n=0,\ldots,P-1$, else False.
			\EndWhile
			\State $y_P=y_{P-1,M}$
			\State \Return $y_P$
			\EndFunction	
		\end{algorithmic}
		\caption{$P$ steps of the monodomain equation in parallel-in-time.}
		\label{algo:hsdc}
	\end{algorithm}
	
	Remark that the HSDC method is fully defined by the number of levels $L$, the numbers of collocation nodes $M_1,\ldots,M_L$ at each level, and the number of processors $P$. To alleviate the notation, henceforth the number of collocation nodes is noted $M=M_1,\ldots,M_L$. For instance, a two-level HSDC method with $M_1=6$ and $M_2=3$ is simply denoted by $M=6,3$.

\section{Stability analysis}\label{sec:stab}
In this section, we begin by numerically assessing the absolute stability of the HSDC method. Subsequently, we demonstrate that the plain SDC method exhibits instability when applied to the monodomain equation, thereby confirming the necessity for HSDC.

\subsection{Stability of HSDC}
We define the linear scalar test equation
\begin{equation}\label{eq:test}
	y'=\lambda_I y+\lambda_E y+\lambda_e y, \qquad y(0)=1,
\end{equation}
where $\lambda_I,\lambda_E,\lambda_e$ represent $f_I,f_E,f_e$, respectively.
We will verify stability of HSDC after one and also after multiple time steps, due to the parallel-in-time approach. Therefore, we define the stability function as
\begin{equation}
	R_P(\lambda_I,\lambda_E,\lambda_e) = y_P,
\end{equation}
where $y_P$ is the approximation of \cref{eq:test} given by HSDC (\cref{algo:hsdc}) after $P$ steps of size $\Delta t=1$. The number of levels $L$ and the collocation nodes $M$ are left out of notation and specified later. Moreover, if $P>1$ we suppose that steps are computed in parallel, hence with one call to \cref{algo:hsdc}.
The stability domain of the method is defined as
\begin{equation}
	\mathcal{S}_P=\lbrace \lambda_I,\lambda_e,\lambda_E\in \Rb\, : \, |R_P(\lambda_I,\lambda_E,\lambda_e)|\leq 1\rbrace.
\end{equation}
Obviously, the exact solution satisfies $|y(P)|\leq 1$ for all $ \lambda_I,\lambda_E,\lambda_e\leq 0$. For HSDC, we cannot expect $|y_P|\leq 1$ for all $\lambda_E\leq 0$ due to the explicit step on the $f_E$ term. However, we don't expect any lower bound on $\lambda_I,\lambda_e\leq 0$, since $f_I,f_e$ are integrated implicitly and exponentially, respectively.

\begin{figure}
	\begin{subfigure}[t]{0.34\textwidth}
		\centering
		\includegraphics[scale=\plotimscaleq,clip,trim=0mm 0mm 12mm 0mm]{\currfiledir 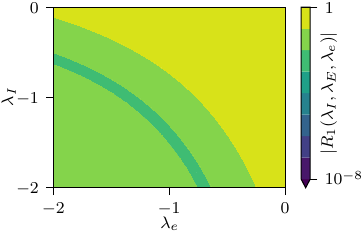}
		\caption{$L=1$, $P=1$, $K=1$.}
		\label{fig:stab_dom_conv_M_6_lo_0_iter_1}
	\end{subfigure}\hfill%
	\begin{subfigure}[t]{0.3\textwidth}
		\centering
		\includegraphics[scale=\plotimscaleq,clip,trim=9mm 0mm 12mm 0mm]{\currfiledir 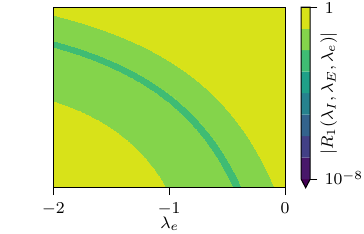}
		\caption{$L=1$, $P=1$, $K=2$.}
		\label{fig:stab_dom_conv_M_6_lo_0_iter_2}
	\end{subfigure}\hfill%
	\begin{subfigure}[t]{0.36\textwidth}
		\centering
		\includegraphics[scale=\plotimscaleq,clip,trim=9mm 0mm 0mm 0mm]{\currfiledir 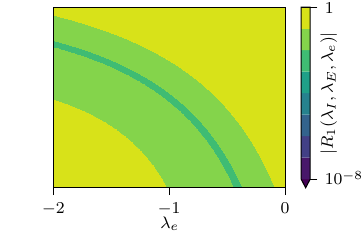}
		\caption{$L=1$, $P=1$, $K=3$.}
		\label{fig:stab_dom_conv_M_6_lo_0_iter_3}
	\end{subfigure}\\
	\begin{subfigure}[t]{0.34\textwidth}
		\centering
		\includegraphics[scale=\plotimscaleq,clip,trim=0mm 0mm 12mm 0mm]{\currfiledir 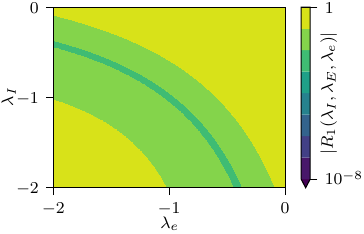}
		\caption{$L=2$, $P=1$, $K=1$.}
		\label{fig:stab_dom_conv_M_63_lo_0_iter_1}
	\end{subfigure}\hfill%
	\begin{subfigure}[t]{0.3\textwidth}
		\centering
		\includegraphics[scale=\plotimscaleq,clip,trim=9mm 0mm 12mm 0mm]{\currfiledir 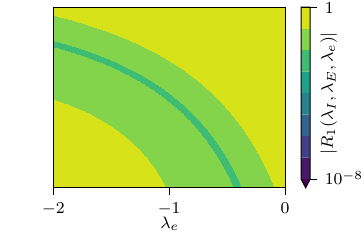}
		\caption{$L=2$, $P=1$, $K=2$.}
		\label{fig:stab_dom_conv_M_63_lo_0_iter_2}
	\end{subfigure}\hfill%
	\begin{subfigure}[t]{0.36\textwidth}
		\centering
		\includegraphics[scale=\plotimscaleq,clip,trim=9mm 0mm 0mm 0mm]{\currfiledir 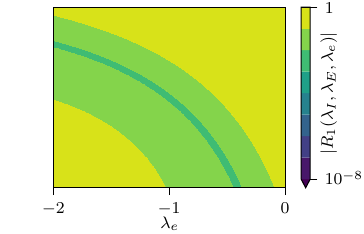}
		\caption{$L=2$, $P=1$, $K=3$.}
		\label{fig:stab_dom_conv_M_63_lo_0_iter_3}
	\end{subfigure}\\
	\begin{subfigure}[t]{0.34\textwidth}
		\centering
		\includegraphics[scale=\plotimscaleq,clip,trim=0mm 0mm 12mm 0mm]{\currfiledir 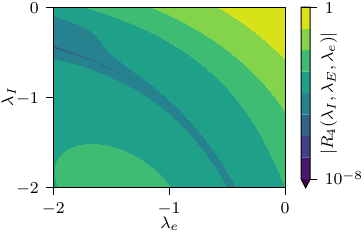}
		\caption{$L=2$, $P=4$, $K=1$.}
		\label{fig:stab_dom_conv_M_63_N_4_lo_0_iter_1}
	\end{subfigure}\hfill%
	\begin{subfigure}[t]{0.3\textwidth}
		\centering
		\includegraphics[scale=\plotimscaleq,clip,trim=9mm 0mm 12mm 0mm]{\currfiledir 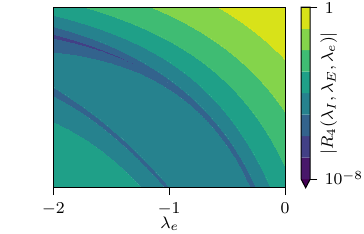}
		\caption{$L=2$, $P=4$, $K=2$.}
		\label{fig:stab_dom_conv_M_63_N_4_lo_0_iter_2}
	\end{subfigure}\hfill%
	\begin{subfigure}[t]{0.36\textwidth}
		\centering
		\includegraphics[scale=\plotimscaleq,clip,trim=9mm 0mm 0mm 0mm]{\currfiledir 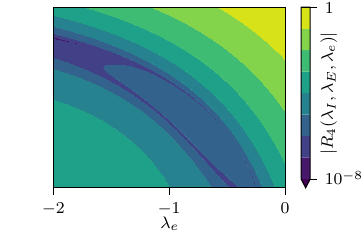}
		\caption{$L=2$, $P=4$, $K=3$.}
		\label{fig:stab_dom_conv_M_63_N_4_lo_0_iter_3}
	\end{subfigure}\\
	\caption{Stability domain for HSDC in a nonstiff regime, with $\lambda_E=0$, $(\lambda_e,\lambda_I)\in [-2,0]\times[-2,0]$ and varying levels $L$, parallel steps $P$, and maximum iterations $K$. The number of collocation nodes is $M=6$ for $L=1$ and $M=6,3$ for $L=2$.}
	\label{fig:stab_dom_nonstiff}
\end{figure}

To verify stability of HSDC, we plot $|R_P(\lambda_I,\lambda_E,\lambda_e)|$ for fixed $\lambda_E\leq 0$ and varying $\lambda_I,\lambda_e\leq 0$.
First, we verify stability in a nonstiff regime, hence we set $\lambda_E=0$ and let $(\lambda_I,\lambda_e)\in [-2,0]\times [-2,0]$.
For a single-level HSDC ($L=1$), with $M=6$ collocation nodes, one time step ($P=1$), fixed relative tolerance $tol=10^{-10}$ and different values of maximum iterations $K=1,2,3$ (cf. \cref{algo:hsdc}) we display the results in \cref{fig:stab_dom_conv_M_6_lo_0_iter_1,fig:stab_dom_conv_M_6_lo_0_iter_2,fig:stab_dom_conv_M_6_lo_0_iter_3}.  
Results for a two-level HSDC with $M=6,3$ are displayed in \cref{fig:stab_dom_conv_M_63_lo_0_iter_1,fig:stab_dom_conv_M_63_lo_0_iter_2,fig:stab_dom_conv_M_63_lo_0_iter_3}. We consider as well a two-level HSDC with $P=4$ parallel steps, results are shown in \cref{fig:stab_dom_conv_M_63_N_4_lo_0_iter_1,fig:stab_dom_conv_M_63_N_4_lo_0_iter_2,fig:stab_dom_conv_M_63_N_4_lo_0_iter_3}. 
We see as in all cases $|R_P(\lambda_I,\lambda_E,\lambda_e)|\leq 1$ and the methods are therefore stable. Furthermore, when $P=1$ (\cref{fig:stab_dom_conv_M_6_lo_0_iter_1,fig:stab_dom_conv_M_6_lo_0_iter_2,fig:stab_dom_conv_M_6_lo_0_iter_3,fig:stab_dom_conv_M_63_lo_0_iter_1,fig:stab_dom_conv_M_63_lo_0_iter_2,fig:stab_dom_conv_M_63_lo_0_iter_3}) the method converges in very few iterations to the same stability domain, indeed for such small values of $\lambda_I,\lambda_E,\lambda_e$ we are in a convergence regime. For $P=4$ parallel steps (\cref{fig:stab_dom_conv_M_63_N_4_lo_0_iter_1,fig:stab_dom_conv_M_63_N_4_lo_0_iter_2,fig:stab_dom_conv_M_63_N_4_lo_0_iter_3}) the scheme converges in few iterations only for those $\lambda_I,\lambda_e$ in a neighborhood of $(0,0)$, but takes longer to converge otherwise. Nevertheless, it is stable for all checked values of $\lambda_I,\lambda_e$.

\begin{figure}
	\begin{subfigure}[t]{0.34\textwidth}
		\centering
		\includegraphics[scale=\plotimscaleq,clip,trim=0mm 0mm 12mm 0mm]{\currfiledir 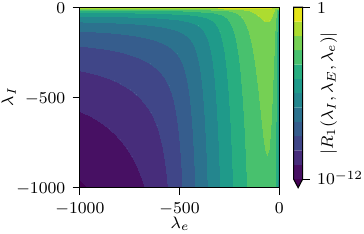}
		\caption{$L=1$, $P=1$, $K=1$.}
		\label{fig:stab_dom_conv_M_6_N_1_lo_-2_iter_1}
	\end{subfigure}\hfill%
	\begin{subfigure}[t]{0.3\textwidth}
		\centering
		\includegraphics[scale=\plotimscaleq,clip,trim=9mm 0mm 12mm 0mm]{\currfiledir 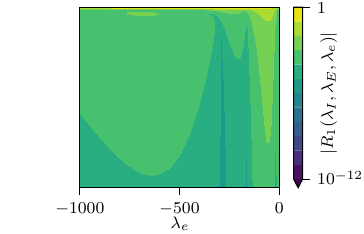}
		\caption{$L=1$, $P=1$, $K=2$.}
		\label{fig:stab_dom_conv_M_6_N_1_lo_-2_iter_2}
	\end{subfigure}\hfill%
	\begin{subfigure}[t]{0.36\textwidth}
		\centering
		\includegraphics[scale=\plotimscaleq,clip,trim=9mm 0mm 0mm 0mm]{\currfiledir 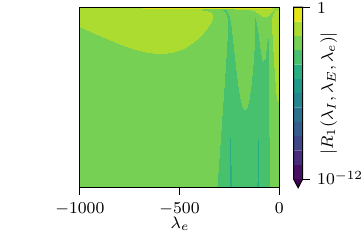}
		\caption{$L=1$, $P=1$, $K=100$.}
		\label{fig:stab_dom_conv_M_6_N_1_lo_-2_iter_100}
	\end{subfigure}\\
	\begin{subfigure}[t]{0.34\textwidth}
		\centering
		\includegraphics[scale=\plotimscaleq,clip,trim=0mm 0mm 12mm 0mm]{\currfiledir 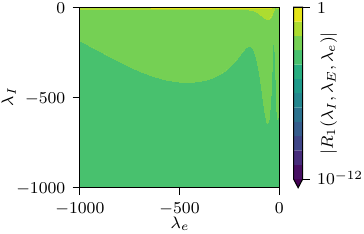}
		\caption{$L=2$, $P=1$, $K=1$.}
		\label{fig:stab_dom_conv_M_63_N_1_lo_-2_iter_1}
	\end{subfigure}\hfill%
	\begin{subfigure}[t]{0.3\textwidth}
		\centering
		\includegraphics[scale=\plotimscaleq,clip,trim=9mm 0mm 12mm 0mm]{\currfiledir 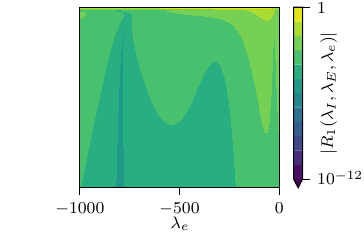}
		\caption{$L=2$, $P=1$, $K=2$.}
		\label{fig:stab_dom_conv_M_63_N_1_lo_-2_iter_2}
	\end{subfigure}\hfill%
	\begin{subfigure}[t]{0.36\textwidth}
		\centering
		\includegraphics[scale=\plotimscaleq,clip,trim=9mm 0mm 0mm 0mm]{\currfiledir 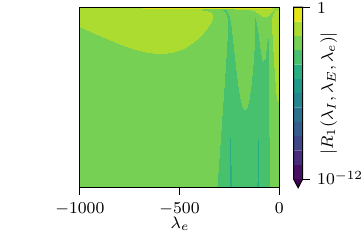}
		\caption{$L=2$, $P=1$, $K=100$.}
		\label{fig:stab_dom_conv_M_63_N_1_lo_-2_iter_100}
	\end{subfigure}\\
	\begin{subfigure}[t]{0.34\textwidth}
		\centering
		\includegraphics[scale=\plotimscaleq,clip,trim=0mm 0mm 12mm 0mm]{\currfiledir 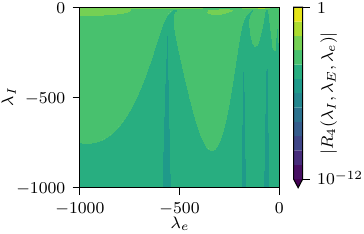}
		\caption{$L=2$, $P=4$, $K=1$.}
		\label{fig:stab_dom_conv_M_63_N_4_lo_-2_iter_1}
	\end{subfigure}\hfill%
	\begin{subfigure}[t]{0.3\textwidth}
		\centering
		\includegraphics[scale=\plotimscaleq,clip,trim=9mm 0mm 12mm 0mm]{\currfiledir 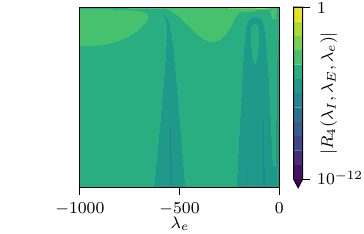}
		\caption{$L=2$, $P=4$, $K=2$.}
		\label{fig:stab_dom_conv_M_63_N_4_lo_-2_iter_2}
	\end{subfigure}\hfill%
	\begin{subfigure}[t]{0.36\textwidth}
		\centering
		\includegraphics[scale=\plotimscaleq,clip,trim=9mm 0mm 0mm 0mm]{\currfiledir 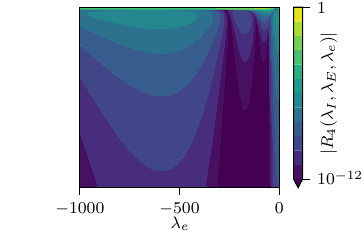}		
		\caption{$L=2$, $P=4$, $K=100$.}
		\label{fig:stab_dom_conv_M_63_N_4_lo_-2_iter_100}
	\end{subfigure}\\
	\caption{Stability domain for HSDC in a stiff regime, with $\lambda_E=-2$, $(\lambda_e,\lambda_I)\in [-1000,0]\times[-1000,0]$ and varying levels $L$, parallel steps $P$, and maximum iterations $K$. The number of collocation nodes is $M=6$ for $L=1$ and $M=6,3$ for $L=2$.}
	\label{fig:stab_dom_stiff}
\end{figure}

In \cref{fig:stab_dom_stiff} we perform similar experiments but in a stiff regime. Hence we fix $\lambda_E=-2$ and let $(\lambda_I,\lambda_e)\in [-1000,0]\times [-1000,0]$. We consider the same combinations of levels $L$, collocation nodes $M$, and processors $P$. However, for the maximum number of iterations we consider $K=1,2,100$, hence in the last case we let the algorithm iterate until convergence.
We observe that at convergence ($K=100$) the serial methods ($P=1$) have a very similar stability domain, compare \cref{fig:stab_dom_conv_M_6_N_1_lo_-2_iter_100,fig:stab_dom_conv_M_63_N_1_lo_-2_iter_100}. Interestingly, the stability domain of the parallel method at convergence, hence \cref{fig:stab_dom_conv_M_63_N_4_lo_-2_iter_100}, has a pattern similar to the serial methods but with much smaller absolute values.
Finally, note that the HSDC scheme is again stable in all considered cases. 

\subsection{Instabilities of SDC for the monodomain equation}
We conclude this section by displaying the stability domain of a naive PinT method for the monodomain equation. Such method would be obtained by using the standard SDC method in combination with a preconditioner based on the popular Rush--Larsen method for the monodomain equation \cite{RuL78}, hence the one defined in \cref{eq:def_P}. Moreover, recall that the plain SDC method is obtained by replacing $\Ln$ in \cref{sec:erk} with $0$.

\begin{figure}
	\begin{subfigure}[t]{0.34\textwidth}
		\centering
		\includegraphics[scale=\plotimscaleq,clip,trim=0mm 0mm 12mm 0mm]{\currfiledir 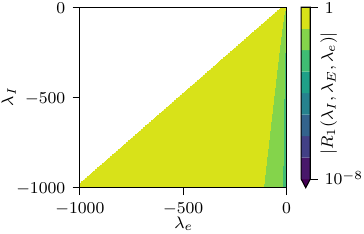}
		\caption{$L=1$, $P=1$, $K=1$.}
		\label{fig:stab_dom_conv_M_6_N_1_lo_-2_iter_1_nonexp}
	\end{subfigure}\hfill%
	\begin{subfigure}[t]{0.3\textwidth}
		\centering
		\includegraphics[scale=\plotimscaleq,clip,trim=9mm 0mm 12mm 0mm]{\currfiledir 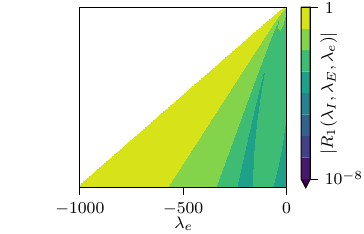}
		\caption{$L=2$, $P=1$, $K=1$.}
		\label{fig:stab_dom_conv_M_63_N_1_lo_-2_iter_1_nonexp}
	\end{subfigure}\hfill%
	\begin{subfigure}[t]{0.36\textwidth}
		\centering
		\includegraphics[scale=\plotimscaleq,clip,trim=9mm 0mm 0mm 0mm]{\currfiledir 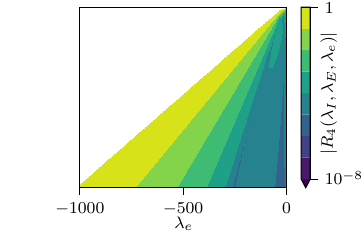}
		\caption{$L=2$, $P=4$, $K=1$.}
		\label{fig:stab_dom_conv_M_63_N_4_lo_-2_iter_1_nonexp}
	\end{subfigure}\hfill%
	\caption{Stability domain of a naive SDC method for the monodomain equation, with $\lambda_E=-2$, $(\lambda_e,\lambda_I)\in [-1000,0]\times[-1000,0]$ and varying levels $L$ and parallel steps $P$. We fix $K=1$ and the number of collocation nodes is $M=6$ for $L=1$ and $M=6,3$ for $L=2$. White regions represent values above $1$, hence instability.}
	\label{fig:stab_dom_stiff_nonexp}
\end{figure}

We show the stability domain of this SDC scheme in \cref{fig:stab_dom_stiff_nonexp}, where we see that already after one iteration the method becomes unstable for $\lambda_e<\lambda_I$. Therefore, the employment of exponential Runge--Kutta methods, as in HSDC, is crucial for stability.

\section{Numerical experiments}\label{sec:numexp}
In this section the accuracy, stability and robustness properties of the HSDC method are assessed numerically. In particular, we first verify the convergence properties of HSDC for single- and multi-level cases and as well in parallel. Second, we investigate the impact of the ionic model and its stiffness on the number of iterations taken by HSDC. Third, we investigate stability and iteration number of HSDC in a realistic setting, considering different step sizes, levels and processors. Finally, we display in details the evolution of the residuals and iterations over time.

\paragraph*{Model setup.} Throughout the next experiments we will solve the monodomain equation with an isotropic conductivity $\sigma=\sigma_i\sigma_e/(\sigma_i+\sigma_e)\,\si{\milli\siemens\per\milli\metre}$, where $\sigma_i=\SI{0.17}{\milli\siemens\per\milli\metre}$ and $\sigma_e=\SI{0.62}{\milli\siemens\per\milli\metre}$ are the intracellular and extracellular electric conductivities. The cell membrane electrical capacitance is set to $C_m=\SI{0.01}{\micro\farad\per\milli\metre\squared}$ and the surface area-to-volume ratio is $\chi=\SI{140}{\per\milli\metre}$. These are the same values as in \cite{NKB11}.
The physical domain will be $\Omega=[\SI{0}{\milli\metre},\SI{100}{\milli\metre}]^d$. The dimension $d=1,2$ and the end time $T$ are specified in the experiments below.

\begin{figure}
	\begin{subfigure}[t]{0.34\textwidth}
		\centering
		\includegraphics[scale=\plotimscaleq,clip,trim=0mm 0mm 4mm 0mm]{\currfiledir 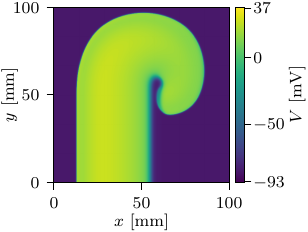}
		\caption{TTP, $t=\SI{0}{\milli\second}$.}
		\label{fig:init_val_TTP}
	\end{subfigure}\hfill%
	\begin{subfigure}[t]{0.30\textwidth}
		\centering
		\includegraphics[scale=\plotimscaleq,clip,trim=8mm 0mm 3.5mm 0mm]{\currfiledir 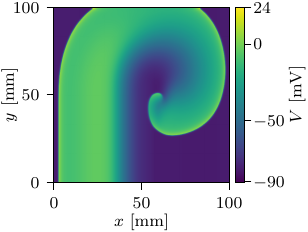}
		\caption{CRN, $t=\SI{0}{\milli\second}$.}
		\label{fig:init_val_CRN}
	\end{subfigure}\hfill%
	\begin{subfigure}[t]{0.33\textwidth}
		\centering
		\includegraphics[scale=\plotimscaleq,clip,trim=8mm 0mm 0mm 0mm]{\currfiledir 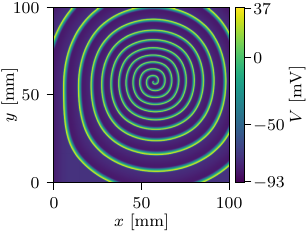}
		\caption{HH, $t=\SI{0}{\milli\second}$.}
		\label{fig:init_val_HH}
	\end{subfigure}\hfill%
	\begin{subfigure}[t]{0.34\textwidth}
		\centering
		\includegraphics[scale=\plotimscaleq,clip,trim=0mm 0mm 4mm 0mm]{\currfiledir 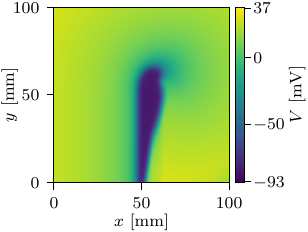}
		\caption{TTP, $t=\SI{100}{\milli\second}$.}
		\label{fig:sol_TTP}
	\end{subfigure}\hfill%
	\begin{subfigure}[t]{0.30\textwidth}
		\centering
		\includegraphics[scale=\plotimscaleq,clip,trim=8mm 0mm 3.5mm 0mm]{\currfiledir 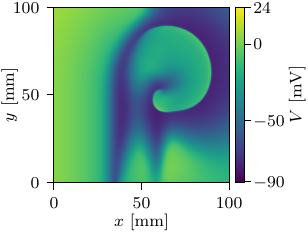}
		\caption{CRN, $t=\SI{100}{\milli\second}$.}
		\label{fig:sol_CRN}
	\end{subfigure}\hfill%
	\begin{subfigure}[t]{0.33\textwidth}
		\centering
		\includegraphics[scale=\plotimscaleq,clip,trim=8mm 0mm 0mm 0mm]{\currfiledir 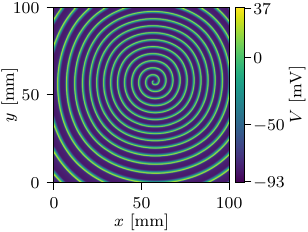}
		\caption{HH, $t=\SI{100}{\milli\second}$.}
		\label{fig:sol_HH}
	\end{subfigure}\hfill%
	\caption{Solutions for different ionic models.}
	\label{fig:sols}
\end{figure}

In most of the forthcoming simulations, we employ the ten-Tusscher--Panfilov (TTP) ionic model \cite{tenPan06}. It is a very stiff model with $m_1=6$ auxiliary variables and $m_2=12$ gating variables.	
In \cref{sec:iter_exp}, where we compare the effect of the ionic models on the convergence of HSDC, we also consider the mildly stiff Courtemance--Ramirez--Nattel (CRN) model \cite{CRN98}, with $m_1=8$ and $m_2=12$, 
and the nonstiff Hodgkin--Huxley (HH) model \cite{HodHux52}, with $m_1=0$ and $m_2=3$. In a typical simulation, with solution away from equilibrium, the spectral radius of the Jacobian of the ionic models' right-hand side is $\rho_{\text{TTP}}\approx 1000$, $\rho_{\text{CRN}}\approx 130$, and $\rho_{\text{HH}}\approx 40$, for TTP, CRN, and HH, respectively. Note that the TTP and CRN are realistic models used in practical applications, while the HH model is mostly used in an academic context. Typical solutions obtained with the three ionic models are displayed in \cref{fig:sols}

For testing the PinT HSDC method we seek for solutions having a fairly constant, in time, ``difficulty''; measured in degrees of stiffness and nonlinearity, for instance. In particular, we want to avoid convergence to an equilibrium, in this way we ensure that all processors have to perform approximately the same amount of work. To do so, we solve monodomain equations showing a self-sustained wave, as depicted in \cref{fig:sols}. There, we also illustrate the initial value, depending on the ionic model, employed in all simulations coming next. Since the initial value already contains a self-sustained wave, there is no need to stimulate the tissue and in all simulations we set $I_{\text{stim}}(t,\bm{x})\equiv 0$.

The domain $\Omega$ is discretized with uniform mesh sizes of common size (cf. \cite{CBC11,NKB11}) $\Delta x_f=100/512\,\si{\milli\metre}\approx \SI{0.2}{\milli\metre}$ or $\Delta x_c=100/256\,\si{\milli\metre}\approx \SI{0.4}{\milli\metre}$ and the Laplacian is approximated with a fourth-order finite difference method.
The time step size $\Delta t$ will usually span a large range of values specified in each experiment.

\paragraph*{Computational setup.} The HSDC method has been implemented in the \texttt{pySDC} \cite{Spe19} library. This library provides an easy to use and flexible framework for the implementation of prototypical codes employing SDC and its derivatives, such as MLSDC, PFASST, and now HSDC. For the particular case of HSDC for the monodomain equation, we implemented the HSDC sweeper, the tailored preconditioner $\mathbf{P}$ and the spatial discretization. 
The linear systems, involving the Laplacian, are solved via discrete cosine transforms (DCT-II). 
The ionic models, due to their high complexity and computational load, are implemented in C++. The code is found in the official \texttt{pySDC} repository \cite{speck_2024_10886512}.
Numerical experiments are performed on the Alps(Eiger) cluster at the Swiss National Supercomputing Centre (CSCS).

\subsection{Convergence order}\label{sec:conv_exp}
In this first experiment we assess the order of convergence of the HSDC scheme, considering different levels, numbers of collocation nodes, and processors.

Recall that a collocation method with $M_1$ Radau IIA collocation nodes converges with order of accuracy $p=2M_1-1$ \cite{HaiWan02}. It is also known that the SDC method with the same nodes converges with order $p=\min(k,2M_1-1)$, where $k$ is the number of iterations \cite{Dutt2000a}.
Instead, an exponential Runge--Kutta collocation method with the same nodes converges with order $p=M_1+1$ (if $M_1>1$) \cite{HocOst2005}. Similarly, ESDC converges with order $p=\min(k,M_1+1)$.
Being HSDC an hybrid method between SDC and ESDC, we expect that it converges with order $p=\min(k,M_1+1)$ (we suppose $M_1>1$, hence $M_1+1\leq 2M_1-1$). However, if the error in the unknowns integrated with SDC is dominant, then orders of accuracy higher than $p=\min(k,M_1+1)$ can as well be observed.

In general, time integration methods for the monodomain equation are based on splitting strategies, for efficiency and stability reasons; therefore, those methods have low order of accuracy due to the splitting error \cite{LinGerJahLoeWeiWie23}. Higher order splitting methods are very scarce in the literature, limited to order $p=4$, and require backward in time steps \cite{CerSpi18}. 
For the membrane equation, hence the monodomain equation without conductivity (\cref{eq:monodomain} with $\sigma=0$), higher order methods are slightly more popular \cite{CouDouPie18,CouLonPie20,DouCouPie18} but also limited to order $p=4$.
In contrast, the HSDC method has theoretically arbitrary order of accuracy. For instance, in this experiment we show that with $M_1=6$ nodes it attains order $p=7$, hence in agreement with the previous considerations. In some cases it converges with order higher than expected, such as $p=9$ with $M_1=6$. Up to the authors knowledge, the maximum convergence order reached so far for the monodomain equation was $p=4$. 

\subsubsection{The serial case}\label{sec:conv_exp_serial}
Here we perform convergence experiments for the multi-level but serial scheme, hence $L\geq 1$ and $P=1$. To do so, we solve the monodomain equation \cref{eq:monodomain} in a two dimensional domain $\Omega=[0,100]\times [0,100]\si{\milli\metre\squared}$, with $T=\SI{1}{\milli\second}$ and the TTP ionic model. We spatially discretize $\Omega$ with the uniform mesh size $\Delta x_c \approx\SI{0.4}{\milli\metre}$.

\begin{figure}
	\begin{subfigure}[t]{\subfigsized\textwidth}
		\centering
		\includegraphics[scale=\plotimscaled, trim=0mm 3mm 0mm 0mm,clip]{\currfiledir 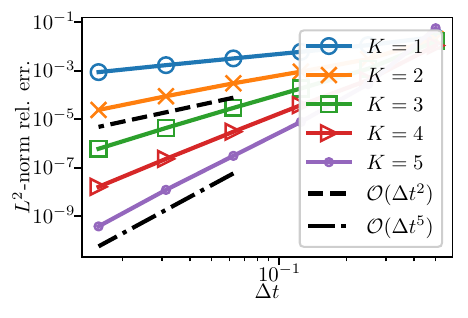}
		\caption{$L=1$, $M=4$.}
		\label{fig:convergence_L_1_M_4}
	\end{subfigure}\hfill%
	\begin{subfigure}[t]{\subfigsized\textwidth}
		\centering
		\includegraphics[scale=\plotimscaled, trim=0mm 3mm 0mm 0mm,clip]{\currfiledir 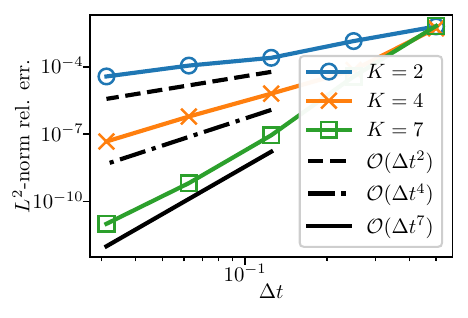}
		\caption{$L=1$, $M=6$.}
		\label{fig:convergence_L_1_M_6}
	\end{subfigure}\hfill%
	\begin{subfigure}[t]{\subfigsized\textwidth}
		\centering
		\includegraphics[scale=\plotimscaled, trim=0mm 3mm 0mm 0mm,clip]{\currfiledir 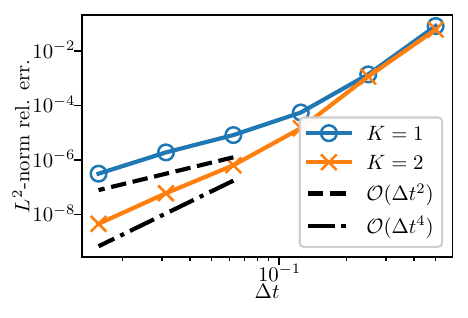}
		\caption{$L=2$, $M=4,2$.}
		\label{fig:convergence_L_2_M_42}
	\end{subfigure}\hfill%
	\begin{subfigure}[t]{\subfigsized\textwidth}
		\centering
		\includegraphics[scale=\plotimscaled, trim=0mm 3mm 0mm 0mm,clip]{\currfiledir 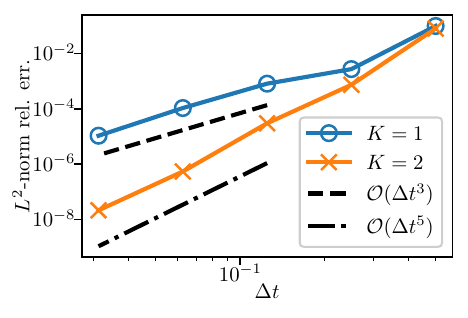}
		\caption{$L=3$, $M=6,3,1$.}
		\label{fig:convergence_L_3_M_631}
	\end{subfigure}\hfill%
	\caption{\cref{sec:conv_exp_serial}. Order of convergence for different levels $L$, collocation nodes $M$ and maximum number $K$ of iterations.}
	\label{fig:convergence_2D}
\end{figure}

To perform the time convergence experiments, we fix the number of levels $L$, of collocation nodes $M$, of maximum iterations $K$, and solve the semi-discrete monodomain equation \cref{eq:sd_sys} for different step sizes $\Delta t$. In \cref{fig:convergence_L_1_M_4} we display the relative errors against a reference solution with respect to $\Delta t$, for a single-level HSDC and $M=4$ collocation nodes, considering different maximum iterations $K$. We see as the order of convergence increases by one with each additional iteration, attaining the theoretically maximal order $p=5$. A similar experiment is performed in \cref{fig:convergence_L_1_M_6}, but with $M=6$ collocation nodes and attaining convergence order $p=7$. In \cref{fig:convergence_L_2_M_42} we show the results for the two-level case, where the order of convergence increases by two at each iteration. Finally, in \cref{fig:convergence_L_3_M_631} the results for the three-level case are displayed, where the order of convergence increases by two or even three at each iteration. Note that for the multi-level cases, we indeed expect to gain more than one order of accuracy per iteration \cite{KreSpe21}.

\subsubsection{The parallel case}\label{sec:conv_exp_parallel}
In this experiment perform convergence experiments for the parallel scheme, hence here $P\geq 1$. 

We consider a one dimensional domain $\Omega=[0,100]\,\si{\milli\metre}$ and set the end time to $T=\SI{16}{\milli\second}$. The domain $\Omega$ is discretized with the uniform mesh size $\Delta x_f\approx \SI{0.2}{\milli\metre}$.

\begin{figure}
	\begin{subfigure}[t]{\subfigsized\textwidth}
		\centering
		\includegraphics[scale=\plotimscaled, trim=0mm 3mm 0mm 0mm,clip]{\currfiledir 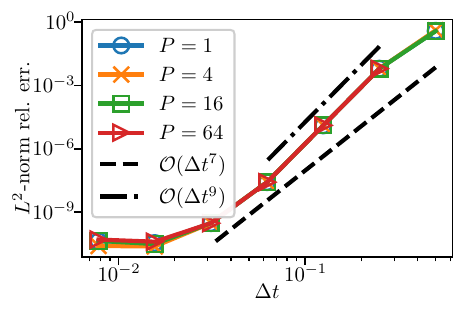}
		\caption{Convergence for different $P$.}
		\label{fig:convergence_L_2_M_63_parallel}
	\end{subfigure}\hfill%
	\begin{subfigure}[t]{\subfigsized\textwidth}
		\centering
		\includegraphics[scale=\plotimscaled, trim=0mm 3mm 0mm 0mm,clip]{\currfiledir 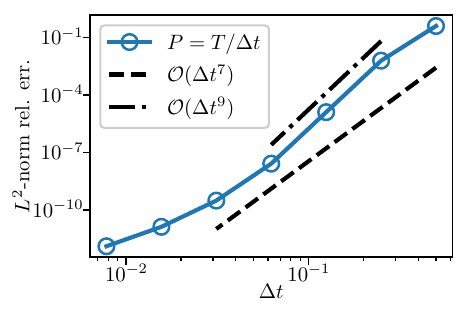}
		\caption{Convergence for $P$ dependent on $\Delta t$.}
		\label{fig:convergence_L_2_M_63_full_parallel}
	\end{subfigure}\hfill%
	\caption{\cref{sec:conv_exp_parallel}. Order of convergence in parallel for $L=2$ and $M=6,3$ collocation nodes. Here HSDC iterates until convergence.}
	\label{fig:convergence_parallel}
\end{figure}

Here we will consider a two-level HSDC with $M=6,3$ collocation nodes. Furthermore, differently from the previous experiment, we do not limit the number of iterations and let the algorithm iterate till convergence, with a relative tolerance $tol=5\cdot 10^{-14}$. 
We display the convergence results for different processors $P=1,4,16,64$ in \cref{fig:convergence_L_2_M_63_parallel}. All cases converge with the expected order of convergence $p=7$ and the accuracy is almost independent from the number of processors. Note that for large step sizes the order of accuracy is even higher than expected ($p=9$ instead of $p=7$). As previously mentioned, orders of accuracies higher than $p=M_1+1$ are possible if the error in the $V_m,w_a$ components of $y$ is dominant; indeed, for these terms the HSDC method corresponds to the SDC method, which is super-convergent and has order $p=2M_1-1$. 

In \cref{fig:convergence_L_2_M_63_full_parallel} we perform a different convergence experiment where the number of processors depends on the step size. More specifically, we set $P=T/\Delta t$, hence all simulation steps are computed simultaneously, disregarding the step size. Since we consider $\Delta t=2^{-k}$, $k=1,\ldots,7$, then the number of processors ranges from $P=32$ to $P=2048$. The same orders of accuracies as in \cref{fig:convergence_L_2_M_63_parallel} are attained.

\subsection{Number of iterations per ionic model}\label{sec:iter_exp}
In this experiment we investigate how the choice of the ionic model impacts the convergence of HSDC. 
More specifically, we solve the monodomain model on a two dimensional domain $\Omega=[0,100]\times [0,100]\si{\milli\metre\squared}$ with a uniform mesh size $\Delta x_f\approx \SI{0.2}{\milli\metre}$ and end time $T=\SI{1}{\milli\second}$.  
Then, for different HSDC settings and ionic models, we compare the average number of iterations (over all time steps) required for convergence, with relative tolerance $tol=5\cdot 10^{-8}$. In this experiment, we fix $P=1$.

\begin{figure}
	\begin{subfigure}[t]{\subfigsized\textwidth}
		\centering
		\includegraphics[scale=\plotimscaled, trim=0mm 3mm 0mm 0mm,clip]{\currfiledir 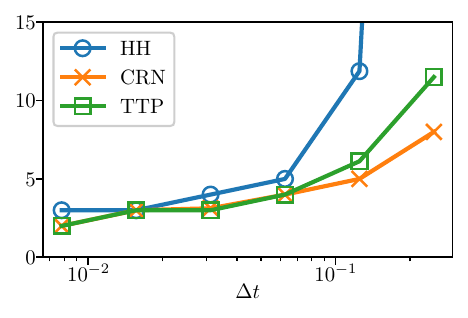}
		\caption{$M=8$ and varying $\Delta t$.}
		\label{fig:iterations_2D_DCT_m_8}
	\end{subfigure}\hfill%
	\begin{subfigure}[t]{\subfigsized\textwidth}
		\centering
		\includegraphics[scale=\plotimscaled, trim=0mm 3mm 0mm 0mm,clip]{\currfiledir 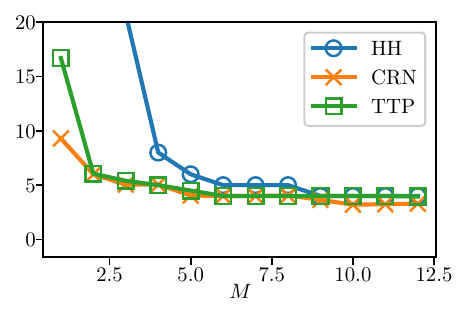}
		\caption{$\Delta t=\SI{0.05}{\milli\second}$ and varying $M$.}
		\label{fig:iterations_2D_DCT_dt_0p05}
	\end{subfigure}\hfill%
	\caption{\cref{sec:iter_exp}. Average number of iterations for different ionic models and varying $\Delta t$ or collocation nodes $M$. Here a single-level HSDC.}
	\label{fig:iterations_2D_one_level}	
\end{figure}

First, in \cref{fig:iterations_2D_DCT_m_8} we consider the three ionic models TTP, CRN, and HH, a single-level HSDC with $M=8$ and vary $\Delta t$. In general, the number of iterations increases with $\Delta t$, as expected, but remains contained for the CRN and TTP models. In contrast, for the HH model there is a significant increase. We believe that this is due to spectral properties: the CRN and TTP models are mostly dissipative, while HH has oscillatory components \cite{Spiteri2010, SpD12}.

Second, in \cref{fig:iterations_2D_DCT_dt_0p05} we fix the step size $\Delta t=\SI{0.05}{\milli\second}$ and let vary the number of collocation nodes $M$ of a single-level HSDC method. Here the average number of iterations increase as $M$ decreases, since the gap between internal stages increases and thus the preconditioner $\bm{P}$ becomes less accurate. Note that again CRN and HH are the less and most demanding models, respectively. Furthermore, for $M>1$ the TTP model is comparable to CRN.

\begin{figure}
	\begin{subfigure}[t]{\subfigsizet\textwidth}
		\centering
		\includegraphics[scale=\plotimscaled, trim=0mm 3mm 0mm 0mm,clip]{\currfiledir 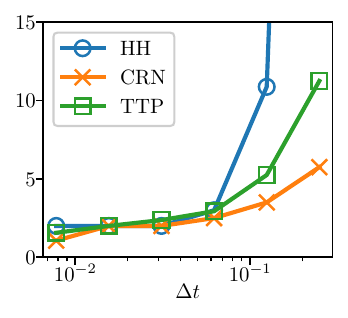}
		\caption{$M=8,6$.}
		\label{fig:iterations_2D_two_levels_M_86}
	\end{subfigure}\hfill%
	\begin{subfigure}[t]{\subfigsizet\textwidth}
		\centering
		\includegraphics[scale=\plotimscaled, trim=0mm 3mm 0mm 0mm,clip]{\currfiledir 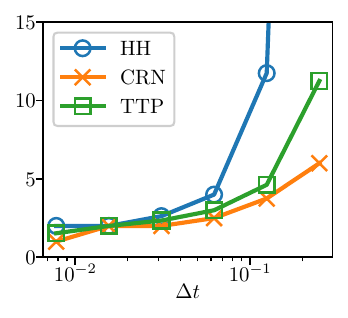}
		\caption{$M=8,4$.}
		\label{fig:iterations_2D_two_levels_M_84}
	\end{subfigure}\hfill%
	\begin{subfigure}[t]{\subfigsizet\textwidth}
		\centering
		\includegraphics[scale=\plotimscaled, trim=0mm 3mm 0mm 0mm,clip]{\currfiledir 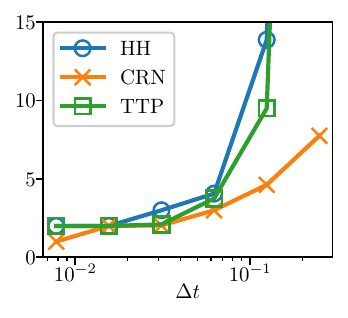}
		\caption{$M=8,2$.}
		\label{fig:iterations_2D_two_levels_M_82}
	\end{subfigure}\hfill%
	\caption{\cref{sec:iter_exp}. Average number of iterations for different ionic models and varying $\Delta t$. Here a two-level HSDC with weak, mild and aggressive coarsening.}	
	\label{fig:iterations_2D_two_levels}
\end{figure}

Finally, we consider a two-level HSDC where we fix the number of fine collocation nodes $M_1$ and consider three choices of coarse collocation nodes $M_2$, leading to weak, mild and aggressive coarsening strategies. We display in \cref{fig:iterations_2D_two_levels_M_86,fig:iterations_2D_two_levels_M_84,fig:iterations_2D_two_levels_M_82} the number of iterations with respect to $\Delta t$ and observe as they tend to increase when the number of coarse collocation nodes $M_2$ decreases, as expected. Furthermore, as in the previous cases, the HH and CRN models are the slowest and fastest to converge, respectively. For step sizes typically used in practice, hence $\Delta t\leq \SI{0.1}{\milli\second}$, the TTP model also takes few iterations to converge.

We conclude pointing out that the slow convergence behavior of HSDC for the HH model in general is not a problem, indeed this model is employed mostly in an academic context and not in practice. In contrast, the robustness of HSDC for the CRN and TTP models is much more relevant, since these are very complex ionic models used in real-life simulations.

\subsection{Parallel-in-time stability}\label{sec:stab_exp}
In \cref{sec:stab}, we showed numerically that for a scalar linear test equation the HSDC method is stable, for both the single- and multi-level case and serial or parallel steps. In this experiment we check stability in a more realistic case, hence solving the monodomain equation \cref{eq:sd_sys} on a two dimensional domain $\Omega=[0,100]\times [0,100]\si{\milli\metre\squared}$ with uniform mesh size $\Delta x_f\approx \SI{0.2}{\milli\metre}$ and the stiff TTP ionic model. 

\begin{figure}
	\begin{center}
		\includegraphics[scale=\plotimscaled, trim=0mm 0mm 0mm 0mm,clip]{\currfiledir 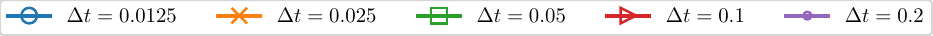}
	\end{center}
	\begin{subfigure}[t]{\subfigsizet\textwidth}
		\centering
		\includegraphics[scale=\plotimscaled, trim=0mm 3mm 0mm 2mm,clip]{\currfiledir 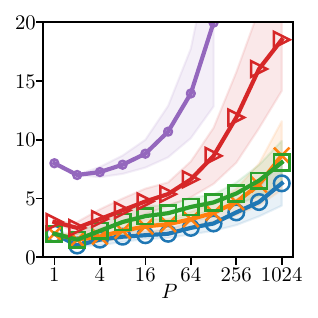}
		\caption{$L=2$, with $M=8,4$.}
		\label{fig:stability_2D_M_84}
	\end{subfigure}\hfill%
	\begin{subfigure}[t]{\subfigsizet\textwidth}
		\centering
		\includegraphics[scale=\plotimscaled, trim=0mm 3mm 0mm 2mm,clip]{\currfiledir 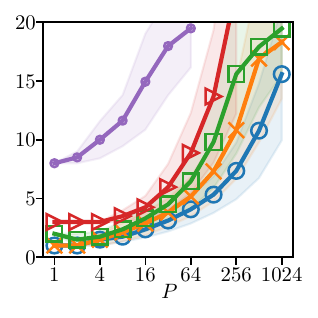}
		\caption{$L=3$, with $M=8,4,2$.}
		\label{fig:stability_2D_M_842}
	\end{subfigure}\hfill%
	\begin{subfigure}[t]{\subfigsizet\textwidth}
		\centering
		\includegraphics[scale=\plotimscaled, trim=0mm 3mm 0mm 2mm,clip]{\currfiledir 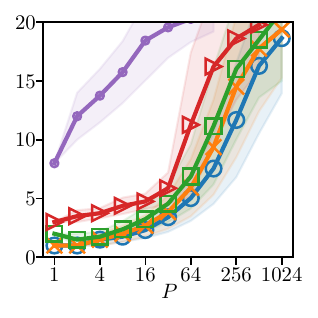}
		\caption{$L=4$, with $M=8,4,2,1$.}
		\label{fig:stability_2D_M_8421}
	\end{subfigure}\hfill%
	\caption{\cref{sec:stab_exp}. Average number of iterations versus number of processors, for different step sizes and levels.}
	\label{fig:stability_2D}
\end{figure}

To carefully asses the stability properties of HSDC, we integrate the problem employing different settings, by varying $\Delta t=0.0125, 0.025,\ldots, 0.2\si{\milli\second}$, number of levels $L=2,3,4$ and parallel steps $P=1,2,\ldots,1024$. For every combination of parameters, we set the end time to $T=P\Delta t$ and let HSDC iterate till convergence. 
The number of iterations required for convergence is an indicator for both efficiency and stability.
In \cref{fig:stability_2D_M_84,fig:stability_2D_M_842,fig:stability_2D_M_8421} we show the results for $L=2,3,4$, respectively. In each figure, we display the average number of iterations (over all time steps) with respect to the number of processors $P$ and the shaded areas indicate standard deviation. Moreover, different curves are associated to a different step size $\Delta t$. 

First, we note as the number of iterations increase with the processors $P$, which is reasonable. For small numbers of processors $P$, the number of iterations decrease with the additional levels, as expected; indeed, in the serial case, the multi-level SDC scheme converges in less iterations than single-level SDC \cite{SpeRupEmmMinBolKra15}. Nevertheless, this property is not always true in the parallel case, indeed for large numbers of processors $P$ iterations increase with the additional levels. The reason comes from the fact that the accuracy of the \textsc{Burn-In} procedure (cf. \cref{algo:hsdc}) depends on the coarsest level. Therefore, if the number of levels is increased by coarsening the coarsest level, like in this experiment, then the \textsc{Burn-In} procedure becomes less accurate and provides rougher initial guesses, leading to an increased number of iterations. This effect is more predominant for long integration intervals, hence many processors.

In \cref{fig:stability_2D_M_84} we see that for $P=2$ there is a decrease in the number of iterations with respect to $P=1$. In general, we expect that iterations increase with the number of processors. However, for $P=1$ we do not call the \textsc{Burn-In} procedure, since it is not needed. In contrast, for $P=2$ it is required and provides initial stage values to both the first and the second processors. This warmer start in the preconditioned Picard scheme leads to a decreased number of iterations. The same behavior is not observed in \cref{fig:stability_2D_M_842,fig:stability_2D_M_8421} since the \textsc{Burn-In} process is not accurate enough to compensate for the increased number of iterations required by the additional processor.

Finally, iterations increase also with the step size $\Delta t$. However, a typical step size for the monodomain equation is $\Delta t\leq \SI{0.05}{\milli\second}$ \cite{CBC11} and here we see that even for $\Delta t\leq \SI{0.1}{\milli\second}$ the number of iterations remains relatively low compared to the corresponding number of processors.

\subsection{Number of iterations and residuals}\label{sec:res_exp}
In this final experiment, we show the number of iterations and the residuals as functions of time. To do so, we discretize again the monodomain equation on the domain $\Omega=[0,100]\times[0,100]\si{\milli\metre\squared}$ with a mesh size $\Delta x_f\approx\SI{0.2}{\milli\metre}$ and a time step size $\Delta t=\SI{0.05}{\milli\second}$.

We solve the semi-discrete equation with HSDC in three different settings. First, a serial and single-level method with $M=8$ collocation nodes. Then, a serial but two-level HSDC with $M=8,4$ collocation nodes. Finally, a parallel two-level HSDC with $M=8,4$, and $P=256$. The end time is set to $T=256\Delta t=\SI{12.8}{\milli\second}$ in all cases.

\begin{figure}
	\begin{center}
		\includegraphics[scale=\plotimscaled]{\currfiledir 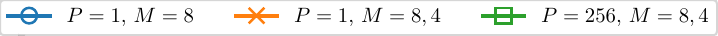}
	\end{center}
	\begin{subfigure}[t]{\textwidth}
		\centering
		\includegraphics[scale=\plotimscaled, trim=0mm 3mm 0mm 2mm,clip]{\currfiledir 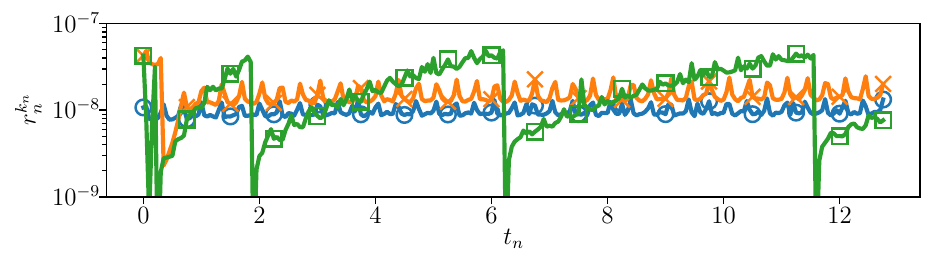}
		\caption{Relative residuals at last iteration $k_n$ over time.}
		\label{fig:residuals_2D_DCT_res_vs_time}
	\end{subfigure}\hfill%
	\begin{subfigure}[t]{\textwidth}
		\centering
		\includegraphics[scale=\plotimscaled, trim=0mm 3mm 0mm 0mm,clip]{\currfiledir 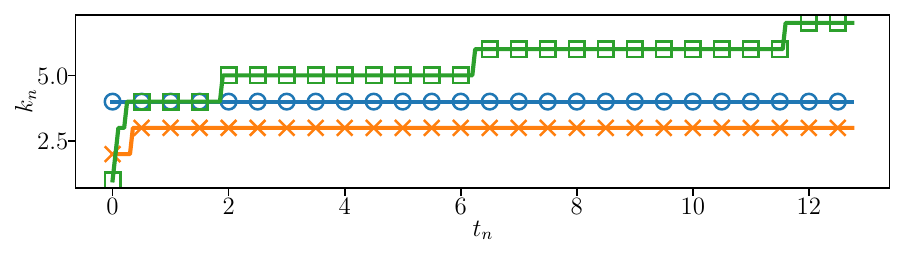}
		\caption{Iterations $k_n$ needed for convergence over time.}
		\label{fig:residuals_2D_DCT_iter_vs_time}
	\end{subfigure}\hfill%
	\begin{subfigure}[t]{\subfigsized\textwidth}
		\centering
		\includegraphics[scale=\plotimscaled, trim=0mm 3mm 0mm 0mm,clip]{\currfiledir 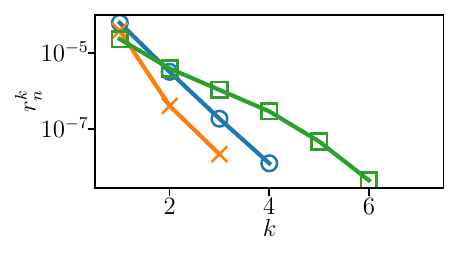}
		\caption{Relative residuals at $t_n=\SI{6.4}{\milli\second}$.}
		\label{fig:residuals_2D_DCT_half_interval}
	\end{subfigure}\hfill%
	\begin{subfigure}[t]{\subfigsized\textwidth}
		\centering
		\includegraphics[scale=\plotimscaled, trim=0mm 3mm 0mm 0mm,clip]{\currfiledir 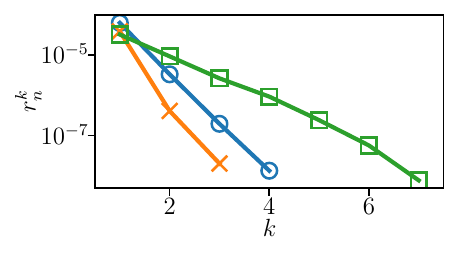}
		\caption{Relative residuals $t_n=\SI{12.8}{\milli\second}$.}
		\label{fig:residuals_2D_DCT_end_interval}
	\end{subfigure}\hfill%
	\caption{\cref{sec:res_exp}. Residuals and iterations count over time.}
	\label{fig:residuals_2D_DCT}	
\end{figure}

We display in \cref{fig:residuals_2D_DCT_res_vs_time} the relative residual 
$r_n^k\coloneqq \Vert\bm{r}_n^{k}\Vert/\Vert\bm{y}_n^{k}\Vert$
with respect to time $t_n$ at iteration $k=k_n$, where $k_n$ is the iterate satisfying the stopping criteria $\Vert\bm{r}_n^{k_n}\Vert<tol\cdot\Vert\by_n^{k_n}\Vert$, i.e. $r_n^{k_n}<tol$, with $tol=5\cdot 10^{-8}$. In \cref{fig:residuals_2D_DCT_iter_vs_time} we display the number of iterations required for convergence, hence $k_n$, with respect to time. Observe as the residual of the single- and two-level serial HSDC methods oscillates around a constant value, with the single-level method having smaller residual but performing more iterations. In contrast, the parallel HSDC method has a residual's pattern that repeatedly grows in time and then suddenly drops when an additional iteration is taken. Note as the parallel method, with $P=256$, resolves the whole time interval $[0,256\Delta t]$ concurrently, but the number of iterations remains comparable to the serial methods and increases only slowly in time.

In \cref{fig:residuals_2D_DCT_half_interval,fig:residuals_2D_DCT_end_interval} we display the relative residual $r_n^k$ with respect to $k$, for fixed time step $n$ at two different times $t_n=128\Delta t=T/2$ and $t_n=256\Delta t=T$. We see as the serial two-level method is the fastest to converge. 

\section{Conclusion}\label{sec:conclusion}
In this study, we developed a novel variant of the spectral deferred correction (SDC) method, named Hybrid Spectral Deferred Correction (HSDC). This method harmoniously blends semi-implicit and exponential SDC techniques, and is specifically tailored to address the complexities of the monodomain equation in cardiac electrophysiology. By leveraging the PFASST framework, HSDC is seamlessly extended to function in a parallel-in-time setting, thereby enhancing computational efficiency by allowing parallelization across the temporal direction. This innovative approach not only preserves the essential stability and robustness needed for such complex simulations but also achieves an arbitrary order of accuracy.

Our comprehensive numerical experiments underscore the efficacy of HSDC, particularly highlighting its capability to adeptly manage challenging ionic models such as the ten-Tusscher--Panfilov. The robustness and accuracy of HSDC are vital for pushing the boundaries of practical real-time simulations in cardiac research.

While the primary focus of this paper has been on the development and validation of the HSDC time discretization scheme, future research should turn towards code optimization and enhancing spatial discretization. This includes integrating space-coarsening strategies within the multi-level framework and exploring the use of the finite element method for spatial discretization. These advancements will undoubtedly enhance the efficiency and broaden the applicability of the method.

\appendix
\section{Matrix coefficients computation}\label{app:comp_aij}
Here we briefly illustrate how to easily compute the matrix coefficients $a_{ij}(z)$ from \cref{eq:def_aij}. The approach proposed here takes inspiration from \cite{Buv20}.	

Given any polynomial $p(s)$ of degree $M-1$, the Fornberg algorithm \cite{For88} can be used to write its derivatives $p^{(k)}(0)$ at $s=0$ in terms of its values at $c_j$, $j=1,\ldots,M$:
\begin{equation}\label{eq:fornberg}
	p^{(k)}(0)=\sum_{l=1}^M w^k_l p(c_{l}).
\end{equation}
From a Taylor expansion at $s=0$ and using \cref{eq:fornberg} and $\ell_j(c_l)=\delta_{jl}$, we have
\begin{equation}\label{eq:fornberg_ell}
	\begin{aligned}
		\ell_j(s)= \sum_{k=1}^{M}\frac{s^{k-1}}{(k-1)!} \ell_j^{(k-1)}(0)
		=\sum_{k=1}^{M}\frac{s^{k-1}}{(k-1)!}\sum_{l=1}^M w^{k-1}_l \ell_j(c_{l})
		=\sum_{k=1}^{M}\frac{s^{k-1}}{(k-1)!}w^{k-1}_j.
	\end{aligned}
\end{equation}
Plugging \cref{eq:fornberg_ell} into \cref{eq:def_aij} yields
\begin{equation}\label{eq:def_aij_computable}
	\begin{aligned}
		a_{ij}(z)= \int_0^{c_i}e^{(c_i-s)z}\ell_j(s)\dif s
		=\sum_{k=1}^M w_j^{k-1}\frac{1}{(k-1)!}\int_0^{c_i}e^{(c_i-s)z}s^{k-1}\dif s
		=\sum_{k=1}^M c_i^k w_j^{k-1}\varphi_k(c_iz)
	\end{aligned}
\end{equation}
where we performed the change of variables $s=c_i r$ and used
\begin{equation}\label{eq:def_phik}
	\varphi_k(z)\coloneqq \frac{1}{(k-1)!}\int_0^1 e^{(1-r)z}r^{k-1}\dif r \quad \text{for }k\geq 1.
\end{equation}
Therefore, the matrix coefficients $a_{ij}(z)$ can be computed using \cref{eq:def_aij_computable}. The coefficients $w_j^k$ are computed with the Fornberg algorithm \cite{For88} and $\varphi_k(z)$ by approximating the integral in \cref{eq:def_phik} with a quadrature rule (we use Gaussian quadrature). The matrix exponential in the integrand of \cref{eq:def_phik} is easily computable since in our case the matrices are diagonal. Notice that computing $a_{ij}(z)$ with formula \cref{eq:def_aij_computable} is more stable than using \cref{eq:def_aij}, as numerical evaluation of the highly oscillatory polynomials $\ell_j(t)$ is avoided.

\bibliographystyle{habbrv}
\bibliography{library}

\end{document}